\documentclass[12pt]{amsart}
\usepackage{amssymb, amsthm, amsmath, url, epsfig, psfrag}
\usepackage{amssymb}
\usepackage[all]{xy}
\usepackage{verbatim}

\def\IZ{{\Bbb Z}}
\def\Z{{\bf Z}}

\def\eproof{$\Box$ \medskip}

\def\PSL{\rm{PSL}}

\def\sideremark#1{\ifvmode\leavevmode\fi\vadjust{\vbox to0pt{\vss
  \hbox to 0pt{\hskip\hsize\hskip1em
  \vbox{\hsize2.5cm\tiny\raggedright\pretolerance10000
  \noindent #1\hfill}\hss}\vbox to8pt{\vfil}\vss}}}

\newtheorem{theorem}{Theorem}[section]
\newtheorem{lemma}[theorem]{Lemma}
\newtheorem{corollary}[theorem]{Corollary}

\newtheorem{prop}[theorem]{Proposition}

\theoremstyle{definition}
\newtheorem{defn}[theorem]{Definition}

\newcommand{\Ht}{\mathbb{H}^3}

\renewcommand{\phi}{\varphi}
\newcommand{\R}{\mathbb{R}}

\newcommand{\rs}{\widehat{\mathbb{C}}}
\newcommand{\C}{\Bbb{C}}

\def\BZ{{\Bbb Z}}

\newcommand{\Out}{{\rm Out}}
\def\Z{\mathbb{Z}}

\begin{document}

\title[Dynamics on $\PSL(2,\C)$-character varieties]{Dynamics on $\PSL(2,\C)$-character varieties: 3-manifolds with toroidal boundary components}
\author{Richard D. Canary and Aaron D. Magid}
\thanks{The first author was partially
supported by NSF grant  DMS - 1006298. The second author was partially supported by the NSF grant DMS-0902764.}
\date\today

\maketitle

\begin{abstract}
Let $M$ be a compact, orientable, hyperbolizable 3-manifold with incompressible
boundary which is not an interval bundle. We study the dynamics of
the action of ${\rm Out}(\pi_1(M))$ on the  relative $\PSL(2,\C)$-character variety $X_T(M)$.
\end{abstract}

\section{Introduction}

We continue the investigation of the action of the outer automorphism group
$\Out(\pi_1(M))$ of the fundamental group of  a compact, orientable, hyperbolizable 
\hbox{3-manifold} $M$ with non-abelian fundamental group on its relative $\PSL(2,\C)$-character variety
$$X_T(M)={\rm Hom}_T(\pi_1(M),\PSL(2,\C)//\PSL(2,C).$$
Here ${\rm Hom}_T(\pi_1(M),\PSL(2,\C))$ is the space of representations 
of $\pi_1(M)$ into $\PSL(2,\C)$  such that  if an element of $\pi_1(M)$ lies in a rank two
free abelian subgroup of $\pi_1(M)$ then its image is either parabolic or the identity.
The set $AH(M)$ of (conjugacy classes of)
discrete, faithful representations is a closed subset of $X_T(M)$ (\cite{Chuckrow, Jorgensen}).  One may think of $AH(M)$ as
the space of marked, hyperbolic 3-manifolds homotopy equivalent to $M$.
The classical deformation theory of Kleinian groups implies that ${\rm Out}(\pi_1(M))$
acts properly discontinuously on the interior ${\rm int}(AH(M))$ of $AH(M)$.
Our main theorem describes a  domain of discontinuity  which is strictly larger than
${\rm int}(AH(M))$ in the case that $M$ has non-empty incompressible
boundary and is not an interval bundle.  Moreover, we characterize when $\Out(\pi_1(M))$ acts
properly discontinuously on an open neighborhood of $AH(M)$.

\begin{theorem}
\label{main}
Let $M$ be a compact, orientable, hyperbolizable $3$-manifold with nonempty incompressible boundary, which is not an
interval bundle. Then there exists an open \hbox{${\rm Out}(\pi_1(M))$}-invariant subset 
$W(M)$ of  $X_T(M)$ such that \hbox{${\rm Out}(\pi_1(M))$}  acts properly discontinuously on $W(M)$, \hbox{${\rm int}(AH(M))$} is a proper subset of $W(M)$, and
$W(M)$ intersects $ \partial AH(M)$.
\end{theorem}

Our proof yields a domain of discontinuity which contains $AH(M)$ in the case where $M$ contains no primitive essential annuli. 
We recall that a properly embedded annulus $A\subset M$ is 
{\em essential} if $\pi_1(A)$ injects into $\pi_1(M))$ and $A$ is not properly homotopic into
$\partial M$. An essential annulus $A$ is {\em primitive} if $\pi_1(A)$ is a
maximal abelian subgroup of $\pi_1(M)$. One may  apply results of Canary-Storm \cite{canary-storm2} to
show that if $M$ does contain a primitive essential annulus, then no such domain of discontinuity can exist.

\begin{corollary}
\label{maincor}
If $M$ is a compact, orientable, hyperbolizable 3-manifold with incompressible boundary 
and non-abelian fundamental group,
then $\Out(\pi_1(M))$ acts properly discontinuously on an open
\hbox{${\rm Out}(\pi_1(M))$}-invariant
neighborhood of $AH(M)$ in $X_T(M)$
if and only if $M$ contains no primitive essential
annuli.
\end{corollary}

If there are no essential annuli with one boundary component contained in  a toroidal boundary
component of $M$, then our main theorem can be extended readily to the full
character variety
$$X(M)={\rm Hom}(\pi_1(M),\PSL(2,\C)//\PSL(2,C).$$

\begin{theorem}
\label{main absolute}
Let $M$ be a compact, orientable, hyperbolizable $3$-manifold with nonempty incompressible boundary, which is not an interval bundle. If $M$ does not contain an essential annulus with
one boundary component contained in a toroidal boundary component of $M$,
then there exists an open \hbox{${\rm Out}(\pi_1(M))$}-invariant subset 
$\hat W(M)$ of  $X(M)$ such that \hbox{${\rm Out}(\pi_1(M))$}  acts properly discontinuously on $\hat W(M)$ and 
$$W(M)=\hat W(M)\cap X_T(M).$$
In particular, $\hat W(M)$ intersects $ \partial AH(M)$.
\end{theorem}

On the other hand, if $M$ does contain an essential annulus with one boundary component
in a toroidal boundary component of $M$, then no point in $AH(M)$ can lie in a domain of
discontinuity for the action of $\Out(\pi_1(M))$.

\begin{prop}
\label{bad annuli}
Let $M$ be a compact, orientable, hyperbolizable  $3$-manifold with nonempty incompressible boundary and non-abelian fundamental group. If $M$ contains an essential annulus with
one boundary component contained in a toroidal boundary component of $M$,
then every point in $AH(M)$ is a limit
of representations in $X(M)$ which are fixed points of infinite order elements of $\Out(\pi_1(M))$.
\end{prop}

Corollary \ref{maincor} then extends to:

\begin{corollary}
\label{maincor absolute}
If  $M$ is a compact, orientable, hyperbolizable  3-manifold with incompressible boundary and non-abelian fundamental group, then  $\Out(\pi_1(M))$ acts properly discontinuously on an open,
\hbox{${\rm Out}(\pi_1(M))$}-invariant
neighborhood of $AH(M)$ in $X(M)$  if and only if $M$  does not contain a primitive essential annulus or
an essential annulus with one boundary component contained in a toroidal boundary
component of $M$.
\end{corollary}

\medskip
\noindent{\bf Historical overview:}
One may view the study of actions of outer automorphism groups on character varieties
as a natural generalization of the study of the action of mapping class groups on
Teichm\"uller spaces.
When $S$ is a closed, oriented, hyperbolic surface, then the mapping class
group ${\rm Mod}(S)$ acts properly discontinuously on Teichm\"uller space
$\mathcal{T}(S)$.
Teichm\"uller space may be identified with a component of the $\PSL(2,\R)$-character
variety $X_2(S)$ of $\pi_1(S)$ and ${\rm Mod}(S)$ is identified with an index two
subgroup of $\Out(\pi_1(S))$. Goldman \cite{goldman-components} showed that $X_2(S)$ has
$4g-3$ components, one of which is identified with $\mathcal{T}(S)$ and another
of which is identified with $\mathcal{T}(\overline{S})$, where $\bar S$ is $S$ given the opposite
orientation. The representations
in the other components are not discrete and faithful.
The mapping class group  preserves
each component of $X_2(S)$ and acts properly discontinuously on the components
associated to $\mathcal{T}(S)$ and $\mathcal{T}(\overline{S})$. Goldman has
conjectured that ${\rm Mod}(F)$ acts ergodically on each of the remaining $4g-5$ components. A resolution
of Goldman's conjecture would give a very satisfying dynamical decomposition.

In this case when $M=S\times [0,1]$ is an untwisted interval bundle,
${\rm Out}(\pi_1(S))$ acts properly discontinuously on the interior $QF(S)$ of $AH(S\times [0,1])$.
$QF(S)$ is the space of quasifuchsian representations, i.e. the convex cocompact
representations.
Goldman  conjectured that the maximal domain of discontinuity  for the
action of $\Out(\pi_1(S))$ on $X(S\times [0,1])$ is exactly $QF(S)$.
For evidence in support of Goldman's conjectures see
Bowditch \cite{bowditch-markoff}, Lee \cite{lee}, Souto-Storm \cite{SS},  and Tan-Wong-Zhang \cite{Tan}.
In particular, it is known that no point in the boundary of $QF(S)$ can lie in a domain of 
discontinuity for ${\rm Out}(\pi_1(S))$ (see Lee \cite{lee}).

Minsky \cite{minsky-primitive} studied the case where $M$ is a handlebody $H_g$ and
exhibited a domain of discontinuity $PS_g$ for $\Out(\pi_1(H_g))=\Out(F_g)$ which is strictly
larger than the the set of convex cocompact representations, and in fact contains
representations which are neither discrete nor faithful. 
(Lee \cite{lee-comp} has generalized Minsky's results in \cite{minsky-primitive} to the more general case where $M$
is a compression body.) Gelander and Minsky
\cite{gelander-minsky} showed that $\Out(F_g)$ acts ergodically on the 
set $R_g$ of  redundant representations. It remains open whether or not $R_g\cup PS_g$ has full
measure in $X(H_g)$.

If $M$ is a twisted interval bundle, Lee \cite{lee} exhibits an explicit domain of
discontinuity $U(M)$ for the action of ${\rm Out}(\pi_1(M))$ on $X(M)$ which is strictly larger than
${\rm int}(AH(M))$ and contains points in $\partial AH(M)$  and points in the complement
of $AH(M)$. Moreover, she shows that if $\rho\in AH(M)-U(M)$, then $\rho$ cannot
lie in a domain of discontinuity for the action of $\Out(\pi_1(M))$ on $X(M)$.

Canary and Storm \cite{canary-storm2} studied the case where $M$
has non-empty incompressible boundary and has no toroidal boundary
components. If $M$ is not an interval bundle,
they again exhibited a domain of discontinuity for the action of $\Out(\pi_1(M))$
which is strictly larger than the interior of $AH(M)$.
Moreover, they showed that there is a domain of discontinuity for the action of
$\Out(\pi_1(M))$ which contains $AH(M)$ if and only if $M$ contains no 
primitive essential annuli.  Our results build on their results.
The major new difficulties in our
case result from the facts that the characteristic submanifold can contain
thickened tori and that
that the mapping class group of $M$ need not have finite
index in ${\rm Out}(\pi_1(M))$ (see Canary-McCullough \cite{CM}).
Our analysis of $\Out(\pi_1(M))$ is necessarily much more intricate than what
was developed in \cite{canary-storm2}.

\bigskip
\noindent
{\bf Outline of argument:}
Our proof relies on exhibiting a finite index subgroup of $\Out(\pi_1(M))$ which is built
from groups of homotopy equivalences associated to components of the characteristic submanifold $\Sigma(M)$ of $M$.
We then consider subgroups of $\pi_1(M)$ which are preserved by these  groups of homotopy equivalences
and study the action of the outer automorphism group of the subgroup of $\pi_1(M)$ on its associated
relative character variety. We can combine these
separate analyses to construct our domain of discontinuity for the action of $\Out(\pi_1(M))$.

\medskip

If $M$ is a compact, orientable, hyperbolizable 3-manifold with incompressible boundary, then
its characteristic submanifold $\Sigma(M)$ consists of solid tori, thickened tori and interval bundles. Johannson
\cite{johannson} showed that every homotopy equivalence can be homotoped so that it preserves $\Sigma(M)$
and restricts to a homeomorphism of $M-\Sigma(M)$. He also showed that only finitely many homotopy
classes of homeomorphisms of $M-\Sigma(M)$ arise, so one can restrict to a finite index subgroup of $\Out(\pi_1(M))$
such that every automorphism is realized by a homeomorphism which restricts to the identity on $M-\Sigma(M)$.
In section \ref{out section}, we build on techniques developed by McCullough \cite{mccullough} and
Canary-McCullough \cite{CM} to construct a finite index subgroup $\Out_0(\pi_1(M))$ of $\Out(\pi_1(M))$
and a short exact sequence 
$$
1 \longrightarrow B \to \Out_0(\pi_1(M)) \xrightarrow{\;\; \Phi \; \; } A \longrightarrow 1
$$
where $A$  is a direct product of mapping class groups of base surfaces of interval bundle components of
$\Sigma(M)$ and cyclic subgroups generated by Dehn twists in vertical annuli in thickened torus components of $M$
and $B$ is the direct product of the free abelian groups generated by Dehn twists in frontier annuli of $\Sigma(M)$
and free abelian groups generated by ``sweeps'' in thickened torus components of $\Sigma(M)$.
Guirardel and Levitt \cite{guirardel-levitt} have recently developed a related short exact sequence for
outer automorphism groups of certain classes of relatively hyperbolic groups.

In section \ref{char annuli} we decompose the frontier of $\Sigma(M)$ into characteristic collections of annuli each of
which is either the entire frontier of a solid torus or thickened torus component of $\Sigma(M)$ or is a single
component of the frontier of an interval bundle component of $\Sigma(M)$. If a characteristic collection of annuli $C$ is the frontier
of a solid torus or a component of the frontier of the interval bundle, we construct a class of free subgroups of $\pi_1(M)$ which
register $C$, in the sense that the group generated by Dehn twists in $C$ injects into the outer automorphism
group of the free subgroup. If $C$ is the frontier of a thickened torus component, then our registering subgroups
are the free product of the fundamental group of the thickened torus component and a free group.
One may use Klein's combination theorem (i.e. the ping-pong lemma), to show that such registering subgroups
exist (see section \ref{existence registering subgroups}).
This entire analysis generalizes the analysis of the mapping
class group ${\rm Mod}(M)$ of $M$ used in Canary-Storm \cite{canary-storm2} in the case that $M$ has no toroidal
boundary components.

In section \ref{smaller deformation spaces}, we show that if $H$ is a registering subgroup, then $\Out(H)$
acts properly discontinuously on the set $GF(H)$  of geometrically finite, minimally parabolic, discrete, faithful representations
and that $GF(H)$ is an open subset of $X_T(H)$. Similarly, if $\Sigma$ is an interval bundle component of
$\Sigma(M)$, we find an open subset
$GF(\Sigma,\partial_1\Sigma)$ of the appropriate relative character variety on which $\Out(\pi_1(M))$ acts
properly discontinuously.
 Our region $W(M)\subset X_T(M)$ is defined to consist of
representations $\rho$ so that for every
characteristic collection of annuli there is a registering subgroup $H$ such that $\rho|_H\in GF(H)$ and for
every interval bundle component $\Sigma$ of $\Sigma(M)$, $\rho|_{\pi_1(\Sigma)}\in GF(\Sigma,\partial_1\Sigma)$.
Proposition \ref{structure of W} establishes that $W(M)$ is an open $\Out(\pi_1(M))$-invariant subset of $X_T(M)$
which contains all discrete, faithful, minimally parabolic representations. In particular, ${\rm int}(AH(M))$ is a proper
subset of $W(M)$ and $W(M) \cap \partial AH(M) \neq \emptyset$.

Proposition \ref{prop disc on W} shows that $\Out(\pi_1(M))$ acts properly discontinuously on $W(M)$, which completes the proof
of our main result. In the proof we consider a sequence $\{\alpha_n\}$ of distinct  elements of $\Out_0(\pi_1(M))$.
We can restrict to a subsequence so that either (1) there exists a thickened torus or interval bundle component  $V$ of $\Sigma(M)$
so that $\alpha_n$ gives rise to a sequence of distinct elements of $\Out(\pi_1(V))$,
or (2) there exists a fixed element $\gamma$ such that $\alpha_n=\beta_n\circ\gamma$ and  a characteristic
collection of annuli $C$ such that  each $\beta_n$ preserves
every registering subgroup $H$ for $C$ and $\{\beta_n\}$ gives rise to a sequence
of distinct elements of $\Out(H)$. Let $D$ be a compact subset of $W(M)$.  We then use the proper discontinuity of the actions of $\Out(\pi_1(V))$ and  $\Out(H)$   to show that $\{\alpha_n(D)\}$ leaves every compact set. 
This allows us to conclude that $\Out_0(\pi_1(M))$ and hence $\Out(\pi_1(M))$
acts properly discontinuously on $W(M)$.

\bigskip
\noindent
{\bf Acknowledgements:} The first author would like to thank Peter Storm for very instructive conversations
which motivated this project.
The second author would like to thank Ian Agol and Neil Strickland for helpful discussions on math overflow.

\section{Background}

\subsection{Deformation spaces of hyperbolic 3-manifolds and character varieties}

Let $M$ be a compact, orientable, hyperbolizable 3-manifold with incompressible boundary. 
Recall that Thurston (see Morgan \cite{Morgan}) proved that  a compact, orientable 3-manifold with non-empty boundary
is hyperbolizable (i.e., 
its interior  admits a complete hyperbolic metric) if and only if it  is atoroidal and irreducible.
We will assume throughout the remainder of the paper
that $M$ has non-abelian fundamental group.
Let $\partial_TM$ denote the non-toroidal boundary components of $M$.

The action of $\Out(\pi_1(M))$ on the character variety $X(M)$ and the relative character
variety $X_T(M)$ is given by 
$$\alpha([\rho])=[\rho\circ \alpha^{-1}]$$
where $\alpha\in \Out(\pi_1(M))$ and $[\rho]\in X(M)$.
(See  Kapovich \cite[Chapter 4.3]{kapovich-book} 
for a discussion of the character variety and the relative character variety.) 

Sitting within $X_T(M)$ is the space $AH(M)$ of (conjugacy classes of) discrete, faithful representations.
If $\rho\in AH(M)$, then $N_\rho=\Ht/\rho(\pi_1(M))$ is a hyperbolic 3-manifold and there exists
a homotopy equivalence $h_\rho:M\to N_\rho$ in the homotopy class determined by $\rho$. One
may thus think of $AH(M)$ as the space of marked hyperbolic 3-manifolds homotopy equivalent to $M$.

The interior ${\rm int}(AH(M))$ of $AH(M)$, as a subset of $X_T(M)$, consists of
discrete, faithful representations which are
geometrically finite and minimally parabolic (see Sullivan \cite{sullivanII}). A representation $\rho\in AH(M)$ is
{\em geometrically finite} if there is a convex, finite-sided fundamental polyhedron for the action of $\rho(\pi_1(M))$ on
$\Ht$. It is {\em minimally parabolic} if $\rho(g)$ is parabolic if and only if $g$ is  a non-trivial element of a
rank two free abelian subgroup of $\pi_1(M)$.

The components of ${\rm int}(AH(M))$ are in one-to-one correspondence with the
set $\mathcal{A}(M)$ of marked, compact, oriented, hyperbolizable 3-manifolds homotopy
equivalent to $M$ and each component is parameterized by an appropriate
Teichm\"uller space (see Bers \cite{bers-survey} or  Canary-McCullough \cite[Chapter 7]{CM} for
complete discussions of this theory).

If $\rho\in AH(M)$, then there is a compact core $M_\rho$ for
$N_\rho$, i.e. a compact submanifold of $N_\rho$ such that the inclusion is a homotopy
equivalence (see Scott \cite{scott}). We may assume that $h_\rho(M)\subset M_\rho$, so that
$h:M\to M_\rho$ is a homotopy equivalence. Formally, $\mathcal{A}(M)$ is the set
of pairs $(M',h')$ where $M'$ is a compact, oriented, hyperbolizable \hbox{3-manifold} and
$h':M\to M'$ is a homotopy equivalence, where $(M_1,h_1)$ is equivalent to
$(M_2,h_2)$ if and only if there is an orientation-preserving homeomorphism
\hbox{$j:M_1\to M_2$} such that $j\circ h_1$ is homotopic to $h_2$. There is
a natural map
$$\Theta:{\rm int}(AH(M))\to \mathcal{A}(M)$$
given by $\Theta(\rho)=[(M_\rho,h_\rho)]$.
Thurston's Geometrization Theorem (see \cite{Morgan}) implies that $\Theta$ is surjective, while Marden's Isomorphism Theorem \cite{marden} implies that the pre-image of any element of
$\mathcal{A}(M)$ is a component. The work of Bers \cite{bers-spaces}, Kra
\cite{kra} and Maskit \cite{maskit-extension}, then
implies that if $(M',h')\in \mathcal{A}(M)$, then
$$\Theta^{-1}(M',h')\cong \mathcal{T}(\partial_T M')$$
where $\mathcal{T}(\partial_T M')$ is the Teichm\"uller space of marked conformal
structures on the non-toroidal components of $\partial M'$.

One may use this parameterization to show that $\Out(\pi_1(M))$ acts
properly discontinuously on ${\rm int}(AH(M))$.

\begin{prop}
\label{prop disc on int}
If $M$ is a compact, orientable, hyperbolizable 3-manifold with incompressible boundary,
then $\Out(\pi_1(M))$ acts properly discontinuously on ${\rm int}(AH(M))$.
\end{prop}

\begin{proof}
Notice that $\Out(\pi_1(M))$ preserves ${\rm int}(AH(M))$. If $Q$ is a component
of ${\rm int}(AH(M))$ and $\Theta(Q)=(M',h')$, let
$${\rm Mod}_+(M',h')\subset \Out(\pi_1(M))$$
denote the set of outer automorphisms which preserve $Q$. An outer automorphism
$\alpha$ lies in \hbox{${\rm Mod}_+(M',h')$} if and only if
\hbox{$(h')_*\circ \alpha\circ (h')_*^{-1}$} is realized by
an orientation-preserving homeomorphism of $M'$.
Thus,  the action of \hbox{${\rm Mod}_+(M',h')$} on \hbox{$Q\cong \mathcal{T}(\partial_T M')$}
 may be identified with the action of a subgroup of ${\rm Mod}(\partial_T M')$ on
 \hbox{$\mathcal{T}(\partial_T M')$}. Therefore, 
since \hbox{${\rm Mod}(\partial_T M')$}
acts properly discontinuously on \hbox{$\mathcal{T}(\partial_TM')$}, ${\rm Mod}_+(M',h')$ acts properly
discontinuously on $Q$. So, $\Out(\pi_1(M))$ acts properly discontinuously on ${\rm int}(AH(M))$.
\end{proof}

\medskip\noindent
{\bf Remark:}  Proposition \ref{prop disc on int} remains true when $M$ has compressible
boundary. The proof above must be altered to take into account that components of
${\rm int}(AH(M))$ are identified with quotients of the relevant Teichm\"uller spaces.

\subsection{The characteristic submanifold} \label{characteristic submanifold background}

If $M$ is a compact, orientable, hyperbolizable 3-manifold with incompressible boundary,
its characteristic submanifold $\Sigma(M)$ contains only
interval bundles, solid tori and thickened tori and the frontier
$Fr(\Sigma(M))$ consists entirely of essential annuli.
The result below recalls the key properties of the characteristic submanifold in our
setting. (The general theory of the characteristic submanifold  was developed by Jaco-Shalen \cite{JS} and Johannson \cite{johannson}). For a discussion of the characteristic submanifold in our hyperbolic setting see
Morgan  \cite[Sec. 11]{Morgan} or Canary-McCullough
\cite[Chap. 5]{CM}).) 

\begin{theorem}
\label{charprop}
Let $M$ be a compact, orientable, hyperbolizable $3$-manifold with incompressible boundary. 
There exists a codimension zero submanifold
\hbox{$\Sigma(M) \subseteq M$} with frontier
\hbox{$Fr(\Sigma(M)) = \overline{\partial \Sigma(M) - \partial M}$} satisfying the following properties: 
\begin{enumerate}
\item Each component $\Sigma_i$ of $\Sigma(M)$ is
either 
\subitem(i)
an interval bundle  over a compact surface with negative
Euler characteristic which intersects $\partial M$ in its associated $\partial I$-bundle,
\subitem(ii)  a thickened torus such that  \hbox{$\partial M\cap \Sigma_i$} contains a 
torus, or
\subitem(iii) a solid torus.
\item The frontier \hbox{$Fr(\Sigma(M))$} is a collection of essential annuli.
\item Any essential annulus in $M$
is properly isotopic  into  $\Sigma(M)$.
\item
If $X$ is a component of \hbox{$M-\Sigma(M)$}, then either $\pi_1(X)$ is non-abelian
or \hbox{$(\overline{X}, Fr(X))\cong (S^1\times [0,1]\times [0,1], S^1\times [0,1]\times \{0,1\}) $}
and $X$ lies between an interval bundle component of
\hbox{$\Sigma(M)$} and a thickened or solid torus component of \hbox{$\Sigma(M)$}.
Moreover, the component of $\Sigma(M)\cup X$ which contains $X$  is not an
interval bundle which intersects $\partial M$ in its associated $\partial I$-bundle.
\end{enumerate}
\noindent A submanifold with these properties is unique up to isotopy, and is called the 
{\em characteristic submanifold} of $M$.
\end{theorem}

\noindent
{\bf Remark:} In Johannson's work, every toroidal boundary component is contained in
some component  of the characteristic submanifold. We use Jaco and Shalen's definition
which requires that no component of the frontier of the characteristic submanifold be
properly homotopic into the boundary. In our setting, one obtains Jaco and Shalen's
characteristic submanifold from Johannson's characteristic submanifold 
by simply removing components which
are regular neighborhoods of toroidal boundary components.

\medskip

Johannson \cite{johannson} proved
that every homotopy equivalence  between compact, orientable, irreducible
3-manifolds with incompressible boundary
may be homotoped so that it
preserves the characteristic submanifold and is a homeomorphism on
its complement.  

\medskip\noindent
{\bf Johannson's Classification Theorem:} (\cite[Theorem 24.2]{johannson})
{\em Let $M$ and $Q$ be compact, orientable, irreducible 3-manifolds with incompressible boundary
and let \hbox{$h:M\to Q$} be a homotopy equivalence. Then $h$ is homotopic to a map
\hbox{$g:M\to Q$} such that
\begin{enumerate}
\item
$g^{-1}(\Sigma(Q))=\Sigma(M)$,
\item
$g|_{\Sigma(M)}:\Sigma(M)\to \Sigma(Q)$ is a homotopy equivalence,
and
\item
$g|_{\overline{M-\Sigma(M)}}:\overline{M-\Sigma(M)}\to \overline{Q-\Sigma(Q)}$ is
a homeomorphism.
\end{enumerate}
Moreover, if $h$ is a homeomorphism, then $g$ is a homeomorphism.
}

\subsection{Ends of hyperbolic 3-manifolds and the Covering Theorem}

In this section, we recall the Covering Theorem which will be used to show that
minimally parabolic, discrete faithful representations lie in our domain of discontinuity
(see Proposition \ref{structure of W}).

We first discuss the ends of the non-cuspidal portion $N^0$ of a hyperbolic 3-manifold with finitely
generated fundamental group. Let $N=\Ht/\Gamma$ be a hyperbolic 3-manifold with finitely
generated fundamental group.
A {\em precisely invariant system of
horoballs} $\mathcal{H}$ for $\Gamma$ is a 
$\Gamma$-invariant collection
of disjoint  open horoballs based at  parabolic fixed points 
of $\Gamma$, such that there is a horoball
based at every parabolic fixed point. It is a consequence of the Margulis
Lemma (see \cite[Theorems D.2.1,  D.2.2]{benedetti-petronio} or \cite[II.E.3, IV.J.17]{maskit-book})
 that every
Kleinian group has a precisely invariant system of horoballs.
Let
$$N^0=(\Ht-\mathcal{H})/\Gamma.$$
Each component of $\partial N^0$ is either an incompressible torus or an incompressible infinite
annulus.
A {\em relative compact core} for $N^0$ is a compact submanifold $R$  of $N^0$ such
that the inclusion of $R$ into $N$ is a homotopy equivalence and $R$ contains
every toroidal component of $\partial N^0$ and intersects every annular component of
$\partial N^0$ in an incompressible annulus.
(McCullough \cite{mccullough} and Kulkarni-Shalen
\cite{kulkarni-shalen} established the existence of a relative compact core).
The {\em ends} of
$N^0$ are in one-to-one correspondence with the components
of $N^0-R$ (see Bonahon \cite[Proposition 1.3]{bonahon}).
An end $E$ of $N^0$ is {\em geometrically finite}
if there exists a neighborhood of $E$ which does not contain any closed geodesics.
Otherwise, the end is called {\em geometrically infinite}.
A hyperbolic 3-manifold $N$ with finitely generated fundamental group
is {\em geometrically finite} if and only if
each end of $N^0$ is geometrically finite (see Bowditch \cite{bowditch-gf} for the equivalence of
the many definitions of geometric finiteness for a hyperbolic 3-manifold).

The Covering Theorem asserts that geometrically infinite ends usually cover finite-to-one.
(The version of the Covering Theorem we state below incorporates the Tameness Theorem of
Agol \cite{agol} and Calegari-Gabai \cite{calegari-gabai}.)

\medskip\noindent
{\bf Covering Theorem:} (Thurston \cite{thurston-notes}, Canary \cite{cover})
{\em
Let $N$ be  a hyperbolic 3-manifold with finitely generated fundamental group  which covers
another hyperbolic 3-manifold $\hat N$ by a local isometry $\pi : N \to \hat N$.
If $E$ is a geometrically infinite end of $N^0$, then either
 
a) $E$ has a neighborhood $U$ such that $\pi$
is finite-to-one on $U$, or
 
b) $\hat N$  has finite volume and has a finite cover $N'$
which fibers over the circle such that, if $N_S$ denotes the
cover of $N'$ associated to the fiber subgroup, then $N$ is
finitely covered by $N_S$.
}

\section{The outer automorphism group} \label{out section}

In this section, we introduce a finite index subgroup $\Out_0(\pi_1(M))$ of $\Out(\pi_1(M))$ and 
show that there exists a short exact sequence
$$1\longrightarrow K(M)\oplus \left(\oplus_i{\rm Sw}(T_i)\right) \longrightarrow \Out_0(\pi_1(M))\longrightarrow \left(\oplus_i D(T_i)\right)\oplus\left(\oplus_j E(\Sigma_j,Fr(\Sigma_j)\right)\longrightarrow 1$$
where $K(M)$ is generated by Dehn twists in annuli of $Fr(\Sigma(M))$,
each ${\rm Sw}(T_i)$ is a free abelian group generated by sweeps supported on
a thickened torus component $T_i$ of $\Sigma(M)$, each $D(T_i)$ is an infinite cyclic
group generated by a Dehn twist in a vertical annulus in a thickened torus component $T_i$ of
$\Sigma(M)$,
and each $E(\Sigma_j,Fr(\Sigma_j))$ is
identified with the mapping class group of the base surface of an interval bundle
component $\Sigma_j$ of $\Sigma(M)$.
Our proof combines work of Johannson \cite{johannson} and Canary-McCullough \cite{CM}
with a new explicit analysis of homotopy equivalences associated to thickened torus
components of the characteristic submanifold.
A similar short exact sequence for a finite index subgroup of the mapping class group 
${\rm Mod}(M)$ was developed by McCullough \cite{VGF} and used in a crucial manner in \cite{canary-storm2}.
Guirardel and Levitt \cite{guirardel-levitt} have developed a related short exact sequence for  finite index subgroups of the outer automorphism groups  of torsion-free, one-ended
relatively hyperbolic groups
which are hyperbolic relative to families of free abelian subgroups.

\subsection{A first short exact sequence and $K(M)$}
Let $\Out_2(\pi_1(M))$ denote the subgroup of $\Out(\pi_1(M))$ consisting of
outer automorphisms which are realized by  homotopy  equivalences which preserves $\Sigma(M)$ and restrict 
to the identity on $M-\Sigma(M)$.
Lemma 10.1.7 (see also Theorem 10.1.9) in \cite{CM} implies that
${\rm Out}_2(\pi_1(M))$ has finite index in ${\rm Out}(\pi_1(M))$.

If $V$ is a component of $\Sigma(M)$, let $E(V,Fr(V))$ be the group of path components 
of the space of homotopy equivalences of $V$ which restrict to homeomorphisms of $Fr(V)$ which are isotopic to the identity. Note that with this definition, a Dehn twist about a frontier annulus of $V$  is a trivial element of  $E(V,Fr(V))$. 
Proposition 10.1.4 in \cite{CM} guarantees that the obvious homomorphism
$$\Psi: {\rm Out}_2(\pi_1(M))\to \oplus E(V_i,Fr(V_i))$$
is well-defined, where the sum is taken over all components of $\Sigma(M)$.
Lemma 10.1.8 in \cite{CM} implies that $\Psi$ is surjective. 

We next show that the kernel $K(M)$ of $\Psi$ is generated by Dehn twists
about frontier annuli. This generalizes Lemma 4.2.2 in McCullough \cite{VGF}. 

\begin{lemma}
\label{K is Dehn}
The kernel $K(M)$ of $\Psi$ is generated by Dehn twists about  the frontier annuli of
$\Sigma(M)$.
\end{lemma}

\begin{proof} If $\alpha$ lies in the kernel of $\Psi$,
then it has a representative which is trivial on $M-\Sigma(M)$, preserves $\Sigma(M)$
and its restriction to $\Sigma(M)$ is homotopic to the identity via a homotopy
preserving $Fr(\Sigma(M))$. We may therefore choose the representative
$h:M\to M$ to
be the identity off of a regular neighborhood $\mathcal{N}$ of $Fr(\Sigma(M))$.

We may choose coordinates so that $\mathcal{N}\cong Fr(\Sigma(M))\times [-1,1]$
and $Fr(\Sigma(M)\subset M$ is identified with $Fr(\Sigma(M)\times \{0\}$ in these
coordinates. We further choose a Euclidean metric on $Fr(\Sigma(M))$ so that
each component is a straight cylinder with geodesic boundary. We can then
homotope $h$ on $\mathcal{N}$ so that each arc of the form $\{ x\}\times [-1,1]$
is taken to a geodesic in the product Euclidean metric on $\mathcal{N}$. It is
easy to check that the resulting map is a product of Dehn twists in the components
of $Fr(\Sigma(M))$.
\end{proof}

We obtain a first approximation to our desired short exact sequence by considering:

$$1 \longrightarrow K(M)\longrightarrow \Out_2(\pi_1(M))\xrightarrow{\;\: \Psi \;\: } \oplus_i E(V_i,Fr(V_i)) \longrightarrow 1.$$

\subsection{The analysis of  $E(V,Fr(V))$} Our next goal is to understand
 $E(V,Fr(V))$ in the various cases. We first recall that $E(V,Fr(V)) $ is finite
when $V$ is a solid torus component of $\Sigma(M)$.

\begin{lemma}
\label{torus E}
{\rm (\cite[Lemma 10.3.2]{CM})}
Let $M$ be a compact, orientable, hyperbolizable 3-manifold with incompressible boundary.
If $V$ is a solid torus component of $\Sigma(M)$, then $E(V,Fr(V))$ is finite.
\end{lemma}

If $\Sigma$ is an interval bundle component of $\Sigma(M)$ with base surface $F$,
then we say $\Sigma$ is {\em tiny} if  its base surface is either
a thrice-punctured sphere or a twice-punctured projective plane. 

The following result combines Propositions 5.2.3 and 10.2.2 in \cite{CM}, see also the
discussion in Section 5 in Canary-Storm \cite{canary-storm2}.

\begin{lemma}
\label{interval bundle E}
Let $M$ be a compact, orientable, hyperbolizable 3-manifold with incompressible boundary.
Suppose $\Sigma$ is an interval bundle component of $\Sigma(M)$ whose base surface $F$ has
negative Euler characteristic.
\begin{enumerate}
\item
$E(\Sigma,\partial \Sigma)$ is identified
with the group ${\rm Mod}_0(F,\partial F)$ of (isotopy classes of)
homeomorphisms of $F$ whose restriction to the boundary is isotopic
to the identity. 
\item $E(\Sigma,Fr(\Sigma))$ injects into
$\Out(\pi_1(\Sigma))$.
\item 
$E(\Sigma,\partial \Sigma)$ is finite if and only if $\Sigma$ is tiny.
\end{enumerate}
\end{lemma}

It remains to analyze the case when $T$ is a thickened torus component of $\Sigma(M)$.
We view $(T,Fr(T))$ as a $S^1$-bundle over $(B,b)$ where $B$ is an annulus and
$b$ is a non-empty collection of arcs in one boundary component $\partial_1 B$ of $B$, so that
$$(T,Fr(T))=(B\times S^1, b\times S^1).$$
 Let $\partial_0 B$ denote the other boundary component. Let
\hbox{$p_1:T\to B$} and \hbox{$p_2:T\to S^1$} be the
projections onto the two factors and let $s:B\to T$ be the section of $p_1$ given by
$s(b)=(b,1)$.
Proposition 10.2.2 in \cite{CM} guarantees
that if \hbox{$f:(T,Fr(T))\to (T,Fr(T))$} is a homotopy equivalence, then it is
homotopic, as a map of pairs, to a fibre-preserving homotopy equivalence
\hbox{$\bar f:(T,Fr(T))\to (T,Fr(T))$}.
Moreover, there is a homomorphism
$$P:E(T,Fr(T))\to E(B,b)$$
given by letting \hbox{$P([f])=[p_1 \circ \bar f\circ s]$}
where $E(B,b)$ is the group of path components of the space
of homotopy equivalences of $B$ which restrict to homeomorphisms of $b$ which are isotopic
to the identity. We will analyze $E(B,b)$ and the kernel of $P$ in order to understand
$E(T,\partial T)$.

\begin{figure}[htbp] \begin{center}
\psfrag{r}[][]{$\partial_1 B$} \psfrag{z}[][]{\footnotesize $\partial_0 B$}
\psfrag{b}[][]{$b_2$}\psfrag{c}[][]{$b_3$}\psfrag{d}[][]{$b_1$}
\psfrag{g}[][]{\small $\gamma_2$}\psfrag{j}[][]{\small $\gamma_3$}\psfrag{k}[][]{\small $\gamma_1$}
\psfrag{a}[][]{\small $\beta$}
\includegraphics[width=2in]{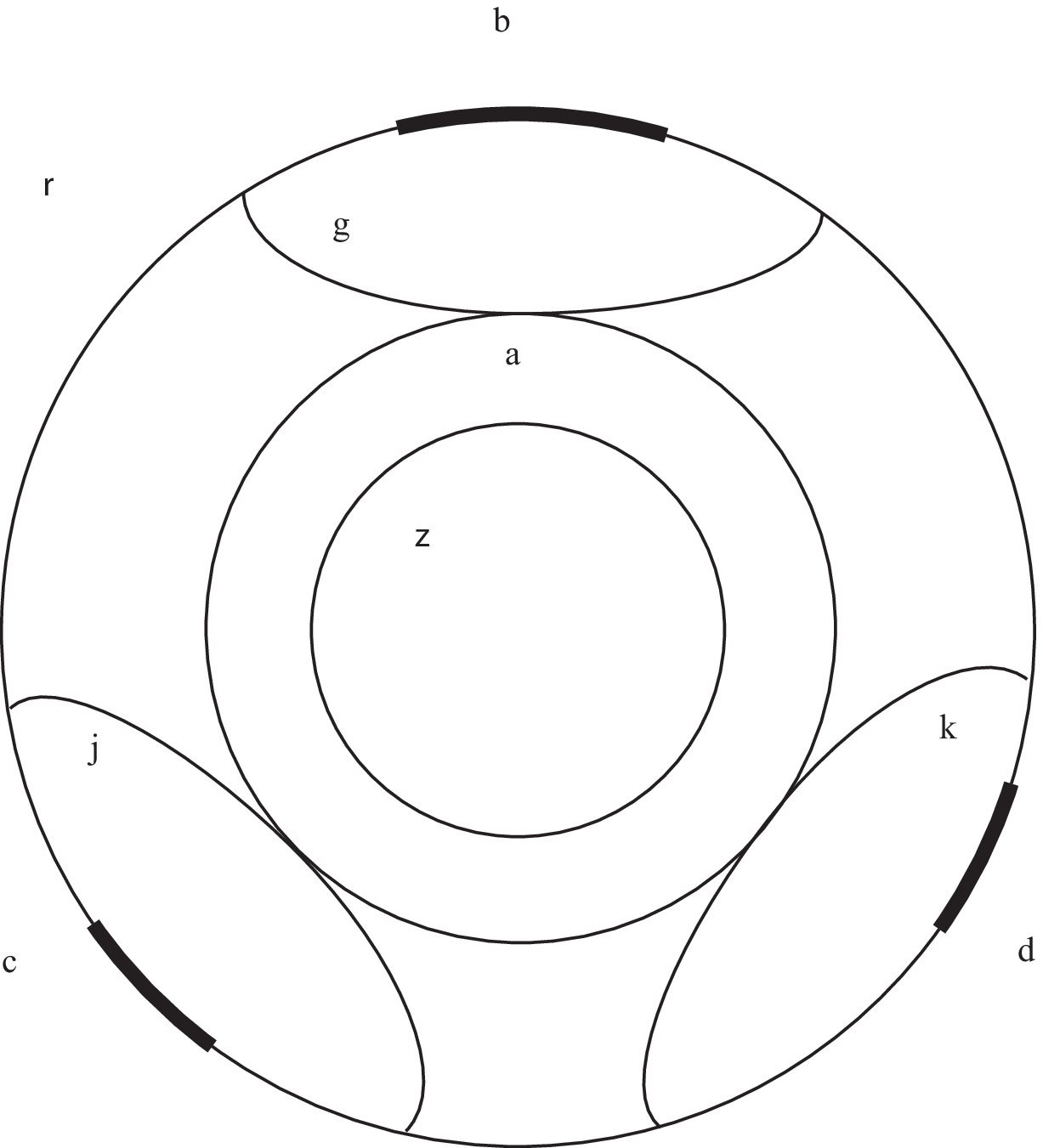}
\hspace{0.9in}
\includegraphics[width=2in]{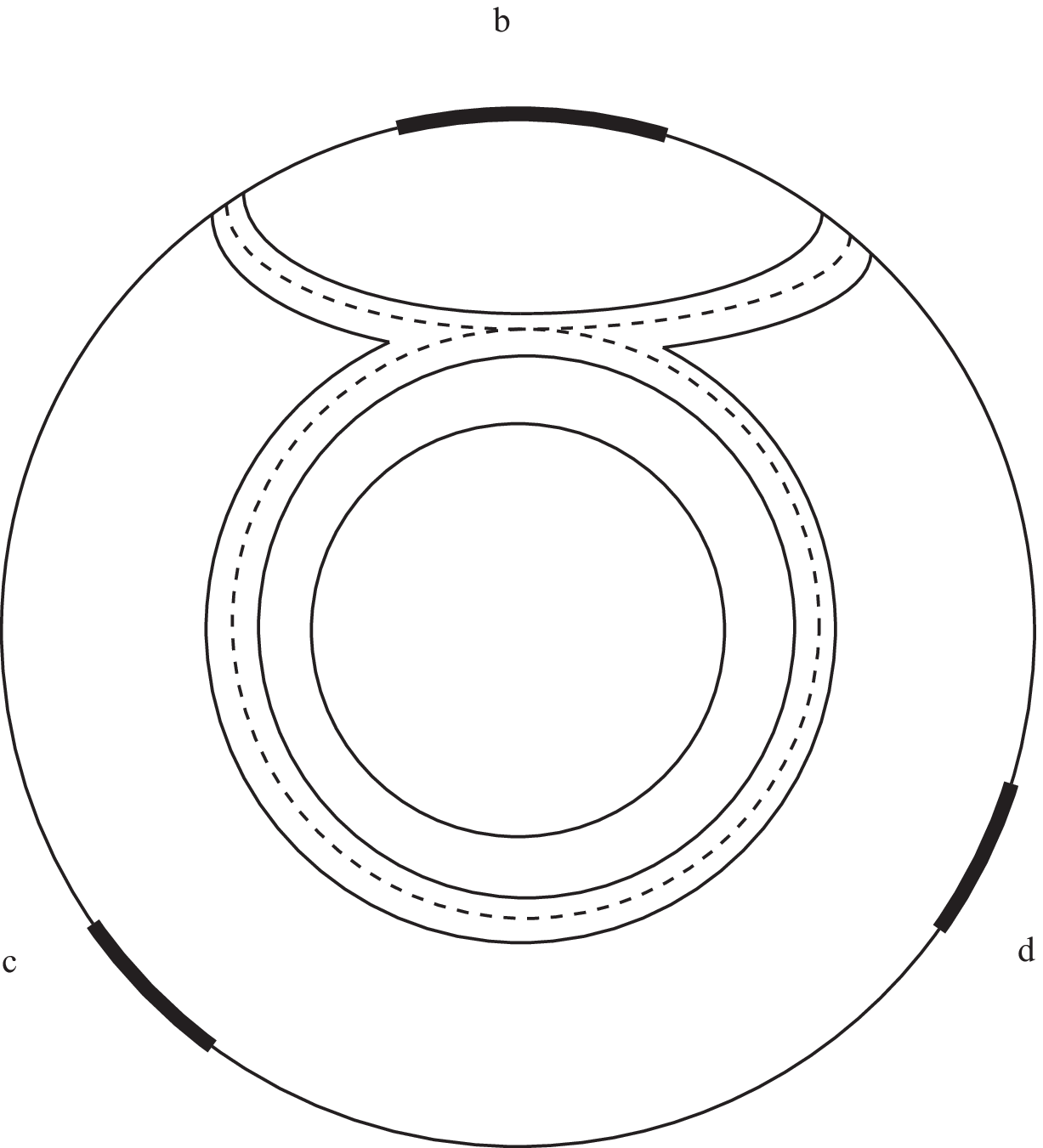}
\caption{
\label{sweep}  The core curve of $B$ is $\beta$ and $b = \{ b_1, b_2, b_3\}$ are three arcs in $\partial_1 B$.   On the right is the image of a neighborhood of $\gamma_2$ after a sweep $h(\gamma_2, \beta)$. 
}
\end{center}\end{figure}

If $\gamma$ is an arc in $B$ with boundary in $\partial_1 B - b$ and $\beta$ is a loop based at a point $x$  on $\gamma$,
then one may define a {\em sweep} $h(\gamma,\beta):(B,b)\to (B,b)$ by requiring that $h$ fixes the complement
of a regular neighborhood $N$ of $\gamma$ and maps a transversal $t$ of $N$ 
through $x$ to \hbox{$t_1*\beta*t_2$} (where $t=t_1*t_2$ and $t_1$ and $t_2$ intersect at $x$).   See Figure \ref{sweep}.
Since $(T,Fr(T))=(B,b)\times S^1$, we may define a sweep
\hbox{$H(\gamma,\beta):(T,Fr(T))\to (T,Fr(T))$} where
\hbox{$H(\gamma,\beta)=h(\gamma,\beta)\times id_{S^1}$}. 
Sweeps are discussed more fully and in greater generality in section 10.2  of \cite{CM}.

Let $b=\{b_1,\ldots,b_n\}$ where each $b_j$ is an arc and $b_j$ is adjacent to $b_{j+1}$
on $\partial B$. Let $\beta$ be a core curve of $B$.
Let $\gamma_j$ be an embedded arc in $B$ joining the components of
$\partial_1 B-b$ adjacent to $b_j$ which intersects $\beta$ at a single point $x_j$.

Let $E_0(B,b)$ denote the index two subgroup of $E(B,b)$ consisting of elements
inducing the identity map on $\pi_1(B)$.

\begin{lemma} \label{generating sweeps}
If $B$ is an annulus and $b=\{b_1,\ldots,b_n\}$ is a non-empty collection of arcs in one component of $\partial B$, then $E_0(B,b)$ is generated by 
$\{ h(\gamma_j,\beta)\}_{j=1}^{n-1}$.
Moreover, $E_0(B,b)\cong \mathbb{Z}^{n-1}$.
\end{lemma} 

We will call $\{ h(\gamma_j,\beta)\}_{j=1}^{n-1}$ a {\em generating system of sweeps} 
for $E_0(B,b)$.

\begin{proof}
We fix an identification of $B$ with $S^1\times [0,1]$.
Choose a collection $\{\lambda_1, \ldots, \lambda_n\}$  of disjoint radial arcs  in $B$ such that each
\hbox{$\lambda_j=\{a_j\}\times [0,1]$} where $a_j\in b_j$. 
See Figure \ref{annulusdiagram}.

\begin{figure}[htbp] \begin{center}
\psfrag{a}[][]{$\partial_1 B$} \psfrag{z}[][]{\footnotesize $\partial_0 B$}
\psfrag{b}[][]{$b_2$}\psfrag{c}[][]{$b_3$}\psfrag{d}[][]{$b_1$}
\psfrag{l}[][]{$\lambda_2$}\psfrag{p}[][]{$\lambda_1$}
\psfrag{n}[][]{$\lambda_3$}
\includegraphics[width=2in]{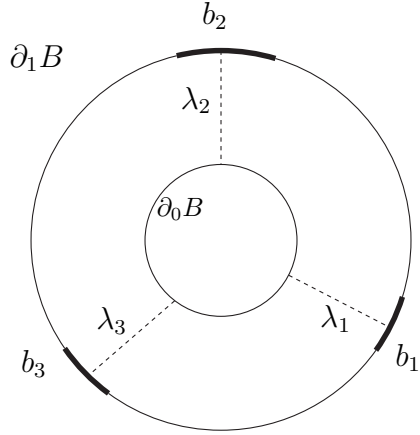}
\caption{
\label{annulusdiagram}  The annulus $(B,b)$ when $n = 3$. The arcs $\lambda_j$ divide $B$ into disks.
}
\end{center}\end{figure}

Given an element of $E_0(B,b)$ we can choose a representative
homotopy equivalence $f$ which is the identity on $b \cup \partial_0 B$. We may further choose $f$ so that $f$
fixes each point on $\lambda_n$. Let $p:S^1\times [0,1]\to S^1$ be projection. Then, for each $j$,
$(p \circ f)_* (\lambda_j)$ will be an element of $\pi_1(S^1,p(a_j))$.
For each $j$, fix an identification of $\pi_1(S^1,p(a_j))$ with $\BZ$.
We define a map  
$$\psi : E_0(B, b) \to \BZ^{n-1}$$
by letting
$$  \psi([f]) =  ((p \circ f)_* (\lambda_1), \ldots, (p \circ f)_* (\lambda_{n-1})).$$

We first show that  $\psi$ is well-defined.  Suppose that $f_1:B \to B$ is another choice of representative of $[f]\in E_0(B,b)$
which is the identity on $b \cup \partial_0 B\cup\lambda_n$.
Since $[f]=[f_1]$ in $E_0(B,b)$,  there is a homotopy $F:B\times [0,1]$ from $f=f_0$ and $f_1$, such
that $F(\cdot,t):\partial_0B\to \partial_0B$ is a homeomorphism isotopic to the identity for all $t$,
We may deform $F$ to a new homotopy, still called $F$, so that $F(\{x\}\times [0,1])$ is a geodesic, in the product
Euclidean metric on $B$, for all $x\in B$. In particular, for all $t$, $F(\cdot,t):\partial_0B\to \partial_0B$ is a rotation.
The fact that $(p \circ f_1)_*(\lambda_n) = (p \circ f_0)_*(\lambda_n)$ implies that $F(\cdot,t):\partial_0B\to \partial_0B$
is the identity for all $t$, so that the homotopy is constant on $\partial_0B$.
Therefore, $(p \circ f_0)_*(\lambda_j) = (p \circ f_1)_*(\lambda_j)$ for $j= 1, \ldots, n-1$, and so $\psi$ is well-defined.

Notice that if $f$ and $g$ are representatives of $[f]$ and $[g]$ in $E_0(B,b)$ which are the identity on
$b \cup \partial_0 B\cup\lambda_n$, then $f\circ g$ is a representative of
$[f][g]$ which is the identity on $b \cup \partial_0 B\cup\lambda_n$, so $\psi$ is a homomorphism.

One may easily check that $\psi(h(\gamma_j,\beta))=(0,\ldots,\pm 1,\ldots 0)$ for all $j$
(where the only non-zero entry is in the $j^{\rm th}$ place. In particular,
$\psi$ is surjective.

The proof will be completed by showing that $\psi$ is injective. 
The collection of arcs $\{ \lambda_j\}$ divides $B$ into $n$ disks $\{D_1,\ldots,D_n\}$. 
Each  $D_i$ has the form $[a_j,a_{j+1}]\times [0,1]$ where indices are taken modulo $n$.
Suppose that $\psi(f)=0$. We may assume as above, that $f$ is the identity on
$b \cup \partial_0 B\cup\lambda_n$.
Since $\psi(f) = 0$, we can further homotope $f$, keeping it the identity on 
$b \cup \partial_0 B\cup\lambda_n$,
so that $f$ fixes $\lambda_j$  for all $j=1,\ldots, n-1$. 
Finally, we  homotope $f$, keeping it the identity on
$b \cup \partial_0 B\cup\lambda_1\cup\cdots\lambda_n$, so that for each $j$ and each $t\in [0,1]$,
$f([a_j,a_{j+1}] ) \times \{t\}$ is a geodesic. This final map must be the identity map, so we have shown
that $[f]=[id]$ in $E_0(B,b)$ which completes the proof.
\end{proof}

If $a$ is an embedded arc in $B$ joining the two boundary components of $B$ such
that $a\cap b=\emptyset$, then $A=p_1^{-1}(a)$ is called a {\em vertical} essential annulus
for $T$. We see that Dehn twists about a vertical annulus generate the kernel of $P$.

\begin{lemma}\label{Dehn twist} Let $T$ be a thickened torus component of $\Sigma(M)$ and let $A$ be a vertical essential annulus in $T$.
Then the kernel of $P$ is generated by the Dehn twist $D_A$ about $A$.
In particular, $ker(P)\cong\BZ$.
\end{lemma}

\begin{proof}
Let $f:T\to T$ be a homotopy equivalence which restricts to a homeomorphism of
$Fr(T)$ which is isotopic to the identity and such that \hbox{$P([f])=[id]$}.
Proposition 10.2.2 of \cite{CM}  allows us to assume
that $f$ is fiber-preserving and that \hbox{$p_1fs=id$}. Then the homotopy class of
$f$ is determined by the homotopy class of the map
\hbox{$p_2fs:B\to S^1$}.  However, every homotopy class of map from $B$ to $S^1$ occurs
when we choose $f=D_A^r$ for some $r$. Therefore, the kernel is generated by $D_A$.
\end{proof}

Let  $E_0(T, Fr(T))=P^{-1}(E_0(B,b))$. Then $E_0(T,Fr(T))$ is an index two subgroup  of $E(T, Fr(T))$ and
there is a short exact sequence 
$$
1 \to ker(P) \to E_0(T, Fr(T)) \to E_0(B,b) \to 1.
$$
Lemma \ref{Dehn twist} shows that $ker(P) \cong \left< D_A \right>  \cong \mathbb{Z}$, while 
Lemma \ref{generating sweeps} shows $E_0(B,b) \cong \mathbb{Z}^{n-1}$ is generated by sweeps
$\{ h(\gamma_j, \beta)\}_{j=1}^{n-1}$.
One may define a section 
$$\sigma:E_0(B,b)\to E_0(T,Fr(T))$$
by setting
$$\sigma(h(\gamma_j,\beta))=H(\gamma_j,\beta)$$
for all $j$. (Notice that  $\sigma$ is  a homomorphism, since any homotopy between maps in $E_0(B,b)$ extends, by simply taking the product
with the identity map on $S^1$, to a homotopy between the corresponding maps in $E_0(T,Fr(T))$.)
Let 
$${\rm Sw}_0(T)=\sigma(E_0(B,b))\cong\Z^{n-1}.$$

For any sweep $H=H(\gamma_j,\beta)$ in $E_0(T,Fr(T))$, $HD_A H^{-1}$ lies in the kernel of $P$.
Since $H$ preserves level sets of $p_2$, we see that
$p_2(HD_AH^{-1})s=p_2(D_A)s$. One then sees that $HD_AH^{-1}=D_A\in E_0(T,Fr(T))$, just as in the proof
of Lemma \ref{Dehn twist}.
Therefore, $E_0(T,Fr(T))$ splits as a direct product
$$E_0(T,Fr(T))=\left< D_A\right> \oplus {\rm Sw}_0(T)\cong \Z\oplus\Z^{n-1}\cong\Z^n.$$

One may combine all of the above analysis to obtain:

\begin{prop}
\label{thickened E}
Let $M$ be a compact, orientable, hyperbolizable 3-manifold with incompressible boundary
and let $T$ be a thickened torus component of $\Sigma(M)$ with base surface $(B,b)$
and $Fr(T)=\{ A_1,\ldots, A_n\}$.
Then  there is 
a subgroup $E_0(T,Fr(T))$ of index  two in $E(T,Fr(T))$ 
which is a rank $n$
free abelian group freely
generated by sweeps $\{H(\gamma_j,\beta)\}_{j=1}^{n-1}$ 
(where $\{ h(\gamma_j,\beta)\}_{j=1}^{n-1}$ is a generating set of sweeps for
$E_0(B,b)$) and a Dehn twist
$D_A$ about a vertical essential annulus $A$.
\end{prop}

\subsection{Assembling the sequence}

We are now ready to define $\Out_0(\pi_1(M))$. If $V$ is a solid torus component
of $\Sigma(M)$, we let
$E_0(V,Fr(V))$ be the trivial group. If $V$ is an interval bundle component of
$\Sigma(M)$, we let $E_0(V,Fr(V))=E(V,Fr(V))$. We have already defined $E_0(V,Fr(V))$
when $V$ is a thickened torus component of $\Sigma(M)$. Then $\oplus_i E_0(V_i,Fr(V_i))$
is a finite index subgroup of $\oplus_i E(V_i,Fr(V_i))$, so
$$\Out_0(\pi_1(M))=\Psi^{-1}\left(\oplus_i E_0(V_i,Fr(V_i))\right)$$
is a finite index subgroup of $\Out_2(\pi_1(M))$ and hence of $\Out(\pi_1(M))$.

If $T$ is a thickened torus component of $\Sigma(M)$, then 
$$E_0(T, Fr(T))=D(T)\oplus {\rm Sw}_0(T)$$
where $D(T)$ is generated by  a Dehn twist about a vertical essential annulus in $T$
and ${\rm Sw}_0(T)$ is generated by sweeps $\{H(\gamma_j,\beta)\}_{j=1}^{n-1}$.
We may extend each $H(\gamma_j,\beta)$ to a homotopy equivalence 
$\hat H(\gamma_j,\beta)$ of $M$ which is the identity on the complement of $T$.
Up to homotopy, we may assume that for all $i \neq j$, the support of $\hat H(\gamma_i,\beta)$ is disjoint from the support of $\hat H(\gamma_j,\beta)$. Also, we may assume the support of $\hat H(\gamma_i,\beta)$ is disjoint from the image of ${\rm supp}( \hat{H}(\gamma_j, \beta))$ under $\hat H(\gamma_j,\beta)$.  It follows that $\hat H(\gamma_i, \beta)$ and $\hat H(\gamma_j, \beta)$ commute. (The argument of Lemma \ref{generating sweeps} can also be adapted to show that
$\hat H(\gamma_i, \beta)$ and $\hat H(\gamma_j, \beta)$ commute.)
One may thus
define a homomorphism
$$s_{T}:{\rm Sw}_0(T)\to {\rm Out}_0(\pi_1(M))$$
by letting $s_T(H(\gamma_j,\beta))=\hat H(\gamma_j,\beta)$ for all $j$.
Since
$\Psi\circ s_{T}$ is the identity map, $s_T$ is an isomorphism onto its image and we define
$${\rm Sw}(T)=s_{T}({\rm Sw}_0(T))\cong \Z^{n-1}.$$ 
We may similarly note that elements of $Sw(T)$ and $K(M)$ commute since one can choose
representatives so that the supports and the images of the supports are disjoint.
(However, we note that $s_{T}$ cannot be extended to a homomorphism defined on all of
$E_0(T,Fr(T))$, since the commutator in $\Out(\pi_1(M))$ of a Dehn twist in a vertical annulus and a sweep is
a Dehn twist in a frontier annulus.)

If we let $\{T_i\}$ be the set of thickened torus components of $\Sigma(M)$ and
$\{\Sigma_j\}$ denote the collection of interval bundle components of $\Sigma(M)$,
then there is an obvious projection map
$$p:\oplus_k E_0(V_k,Fr(V_k))\to  \left(\oplus_i D(T_i)\right)\oplus\left(\oplus_j E(\Sigma_j,Fr(\Sigma_j)\right).$$
Then we consider the map
$$\Phi=p\circ \Psi:\Out_0(\pi_1(M))\to \left(\oplus_i D(T_i)\right)\oplus\left(\oplus_j E(\Sigma_j,Fr(\Sigma_j)\right)$$
which has kernel $K(M)\oplus \left(\oplus_i{\rm Sw}(T_i)\right)$.

We summarize the above discussion in the following proposition.

\begin{prop}
\label{SES}
Let $M$ be a compact, orientable, hyperbolizable 3-manifold with incompressible
boundary. Then there exists a finite index subgroup $\Out_0(\pi_1(M))$ of
$\Out(\pi_1(M))$ and a short exact sequence
$$1\longrightarrow K(M)\oplus \left(\oplus_i{\rm Sw}(T_i)\right) \longrightarrow \Out_0(\pi_1(M)) \xrightarrow{\;\: \Phi \; \:} \left(\oplus_i D(T_i)\right)\oplus\left(\oplus_j E(\Sigma_j,Fr(\Sigma_j)\right)\longrightarrow 1$$
\end{prop}

\section{Characteristic collections of annuli and registering subgroups} \label{char annuli}

We will divide the annuli in $Fr(\Sigma(M))$ into collections,
called characteristic collections of annuli. Each isotopy class of annulus in
$Fr(\Sigma(M))$ will appear in exactly one collection.
A {\em characteristic collection of annuli} $C$ for $\Sigma(M)$ is either
\begin{enumerate}
\item
the collection of frontier annuli in a solid torus component of $\Sigma(M)$,
\item
a component of the frontier of an interval bundle component of $\Sigma(M)$
which is not isotopic into either a solid torus or thickened torus component of $\Sigma(M)$, or
\item
the collection of frontier annuli in a thickened torus component of $\Sigma(M)$.
\end{enumerate}

Let $\{C_1,\ldots,C_m\}$ denote the collections of characteristic annuli for $M$.
Let $K(C_j)$ be the subgroup of $K(M)$ generated by Dehn twists about the annuli in
$C_j$. Notice that
$$K(M)\cong \oplus K(C_j).$$ 

We extend this decomposition of $K(M)$ into a decomposition of ${\rm ker}(\Phi)$.
If $C$ is the collection of frontier annuli of a thickened torus $T$, we define
$\hat K(C)=K(C)\oplus {\rm Sw}(T)$. Otherwise, we define $\hat K(C)=K(C)$. With this
convention,
$${\rm ker}(\Phi)=K(M)\oplus \left(\oplus_i{\rm Sw}(T_i)\right)=\oplus_j \hat K(C_j).$$
If $C$ is a characteristic collection of annuli, then we may define the
projection map
$$q_C:{\rm ker}(\Phi)\to \hat K(C).$$

We next define subgroups of $\pi_1(M)$ which ``register'' the action of $\hat K(C)$ on
$\pi_1(M)$, in the sense that the subgroup is preserved by any element of
$\Out(\pi_1(M))$ and $\hat K(C)$ injects into the outer automorphism group of the subgroup.

Let $M_C=M-\mathcal{N}(C_1 \cup C_2 \cup \ldots \cup C_m)$ be the complement of a
regular neighborhood of the characteristic collections of annuli. If $X$ is a component of $M_C$, then
$X$ is either properly isotopic to  a component of $\Sigma(M)$ or  to 
a component of \hbox{ $M-\Sigma(M)$}. In particular, $\pi_1(X)$ is non-abelian if
it is not properly isotopic to a solid torus or thickened torus component of $\Sigma(M)$.
Moreover, no two adjacent components of $M_C$ have abelian fundamental group.

First suppose that \hbox{$C=Fr(V)=\{A_1,\ldots,A_l\}$} where $V$ is 
a solid torus component of \hbox{$\Sigma(M)$}. For each \hbox{$i=1,\ldots,l$}, let $X_i$ be
the component of $M_C-V$ abutting $A_i$. 
Let $a$ be a core curve for $V$ and let $x_0$ be a point on $a$.
We say that a subgroup  $H$ of \hbox{$\pi_1(M,x_0)$} is  {\em $C$-registering}
if there exist, for each \hbox{$i=1,\ldots,l$}, a loop $g_i$ in \hbox{$T_j\cup X_i$} based at $x_0$ intersecting $A_i$ exactly twice, such that
$$H=<a>*<g_1>*\cdots*<g_l>\cong F_{l+1}$$ 

Now suppose that $C=\{A\}$ is a component of $Fr(\Sigma)$ where $\Sigma$ is an interval 
bundle component of $\Sigma(M)$. Let $a$ be a core curve for $A$ and let $x_0$ be
a point on $a$. We say that a subgroup $H$ of \hbox{$\pi_1(M,x_0)$} is
{\em $C$-registering} if there exist
two loops $g_1$ and $g_2$ based
at $x_0$ each of whose interiors misses $A$, and which lie in the two
distinct components of \hbox{$M_C$} abutting $A$, such that
$$ H=<a>*<g_1>*<g_2>\cong F_{3}$$ 

Finally, suppose that \hbox{$C=Fr(T)=\{A_1,\ldots,A_l\}$} where $T$ is  a thickened torus
component of $\Sigma(M)$.  For each \hbox{$i=1,\ldots,l$}, let $X_i$ be
the component of $M_C-T$ abutting $A_i$. Pick $x_0\in T$.
We say that a subgroup  $H$ of \hbox{$\pi_1(M,x_0)$} is  {\em $C$-registering}
if there exist, for each \hbox{$i=1,\ldots,l$}, a homotopically non-trivial loop $g_i$ in 
\hbox{$T\cup X_i$} based at $x_0$ intersecting $A_i$ exactly twice, such that
$$H=\pi_1(T,x_0)*<g_1>*\cdots*<g_l>\cong (\IZ\oplus\IZ)*F_l$$ 

If $H$ is a subgroup of $\pi_1(M)$, then there is an obvious map 
$$r_H:X(M)\to X(H)={\rm Hom}(H,\PSL(2,\C))//\PSL(2,\C)$$
given by taking $\rho$ to $\rho|_H$.

The following lemma records the key properties of registering subgroups.

\begin{lemma}
\label{characteristic}
Let $M$ be a compact, orientable, hyperbolizable 3-manifold with incompressible
boundary. If $C$ is a characteristic collection of annuli for $M$ and
$H$ is a $C$-registering subgroup of $\pi_1(M)$,
then $H$ is preserved by each element of $\hat K(C)$ and there is a natural
injective homomorphism
$$s_H:\hat K(C)\to {\rm Out}(H).$$
Moreover,  if \hbox{$\eta \in {\rm ker}(\Phi)=K(M)\oplus \left(\oplus_i{\rm Sw}(T_i)\right)$},
then 
$$r_H(\rho\circ \eta)=r_H(\rho)\circ s_H(q_C(\eta))$$
for all \hbox{$\rho\in X(M)$}, where $q_C:{\rm ker}(\Phi)\to\hat K(C)$ is the projection map. 
\end{lemma}

\begin{proof}
If $C$ lies in the frontier of a solid torus or interval bundle component of $\Sigma(M)$, this was
established as Lemma 6.1 in \cite{canary-storm2}.

Now suppose that $C=Fr(T)=\{A_1,\ldots,A_n\}$ where $T$ is a thickened torus component of
$\Sigma(M)$ and that 
$H$ is a \hbox{$C$-registering} subgroup of
$\pi_1(M,x_0)$ (where $x_0\in {\rm int}(T)$) generated by $\pi_1(T,x_0)$ and $\{ g_1,\ldots,g_n\}$.
Then $K(C)\cong \IZ^{n-1}$ and is generated by the Dehn twists $\{D_{A_1},\ldots, D_{A_{n-1}}\}$
and ${\rm Sw}(T)\cong \IZ^{n-1}$ and is generated by sweeps 
$\{\hat H(\gamma_1,\beta),\cdots,\hat H(\gamma_{n-1},\beta)\}$.
We choose generators $a$ and $b$ for $\pi_1(T,x_0)$ so that $a$ is homotopic to the core
curve of $A_1$ and $b$ is the core curve $\beta$ of the annulus $B$ 
(here we adapt the notation of Proposition \ref{thickened E}).
One may check that each $(D_{A_k})_*$ preserves $H$, fixes 
$\pi_1(T,x_0)$ and each $g_i$ where $i\ne k$, 
and maps $g_k$ to $ag_ka^{-1}$.
Similarly, each $(\hat H(\gamma_k,\beta))_*$ preserves $H$ and fixes
$\pi_1(T,x_0)$ and each $g_i$ where $i\ne k$ and takes $g_k$ to
$bg_kb^{-1}$. 

This explicit description of each generator allows one to immediately check that 
$\hat K(C)$ preserves $H$ and injects into $\Out(H)$.
\end{proof}

\section{The existence of registering subgroups} \label{existence registering subgroups}

In this section, we prove that every characteristic collection of annuli admits
a registering subgroup. We will make use of the existence of minimally parabolic
hyperbolic structures on $M$ to do so.

In the case that  the characteristic collection of annuli lies in the frontier of
a solid torus or interval bundle component of $\Sigma(M)$,
the proof of Lemma 8.3 in \cite{canary-storm2} immediately yields:

\begin{lemma}
\label{solid torus registering}
Suppose that $M$ is a compact, orientable, hyperbolizable 3-manifold with incompressible
boundary and $C$ is a characteristic collection of annuli for $M$ such
that  either (1) $C=Fr(V)$ for a solid torus  component $V$ of $\Sigma(M)$ or
(2) $C$ is a component of the frontier of an interval bundle component of $\Sigma(M)$.
If \hbox{$\rho\in AH(M)$} and $\rho(\pi_1(C))$ is purely hyperbolic,
then there exists a $C$-registering subgroup $H$ of
$\pi_1(M)$ such that $\rho|_H$ is  discrete, faithful, geometrically finite
and purely hyperbolic (and therefore Schottky). 
\end{lemma}

We may give a variation on the argument in \cite{canary-storm2} to prove:

\begin{lemma}
\label{thickened torus registering}
Suppose that $M$ is a compact, orientable, hyperbolizable 3-manifold with incompressible
boundary and $C$ is a characteristic collection of annuli for $M$ such
that   $C=Fr(T)$ for a thickened torus  component $T$ of $\Sigma(M)$.
If \hbox{$\rho\in AH(M)$},
then there exists a $C$-registering subgroup $H$ of
$\pi_1(M)$ such that $\rho|_H$ is  discrete, faithful, geometrically finite
and minimally parabolic.
\end{lemma}

\begin{proof}
Let $\rho\in AH(M)$ and let $C = \{ A_1,\ldots, A_l\}$. Let $X_i$ be the component of $M_C$ abutting
$A_i$. Pick $x_0\in T$ and, in order to be precise, let $\pi_1(T)$ denote $\pi_1(T,x_0)\subset \pi_1(M,x_0)$.
Since $\rho\in AH(M)$, $\rho(\pi_1(T))$ consists of parabolic elements fixing a common fixed point
$p\in\rs$.  We may assume that $p=\infty$ and pick a fundamental domain $F$ for  the action of $\rho(\pi_1(T))$ 
on $\C$ which is a quadrilateral.
Since $\pi_1(X_i\cup T,x_0)$ is not abelian, we can find 
$\gamma_i \in \pi_1(X_i\cup T_0,x_0)$ such that $\rho(\gamma_i)$ is a hyperbolic
element with both fixed points contained in
the interior of $F$.  If $i\ne j$, then $\pi_1(X_i\cup T_0,x_0)\cap \pi_1(X_j\cup T_0,x_0)=\pi_1(T,x_0)$,
so $\gamma_i$ and $\gamma_j$ have distinct fixed points. One may then find a collection
$\{D_1^{\pm},\ldots,D_l^\pm\}$ of $2l$ disjoint disks in the interior of $F$ and integers $\{ s_1,\ldots, s_l\}$,
so that, for each $i$, $\gamma_i^{s_i}$ takes the interior of $D_i^-$ homeomorphically onto the exterior of $D_i^+$.
For each $i$, let $g_i$ be a curve in $X_i\cup T$ which intersects $A_i$ exactly twice and  represents
$\rho^{-1}(\gamma_i^{s_i})$. Let
$$H=\left<\pi_1(T),g_1,\ldots,g_l\right>\subset \pi_1(M).$$
Klein's Combination Theorem (see \cite[Theorem A.13, Theorem C.2]{maskit-book}) then implies that 
$\rho$ is geometrically finite and minimally parabolic and that
$$\rho(H)\cong \rho(\pi_1(T)) * \left<\gamma_1^{s_1} \right> *\cdots * \left<\gamma_l^{s_l} \right>.$$
Therefore, $H$ is a registering subgroup with the desired properties.
\end{proof}

Thurston's Hyperbolization Theorem, see Morgan \cite{Morgan}, implies
that there exists a geometrically finite, minimally parabolic element $\rho\in AH(M)$.  
If $C$ is a characteristic collection of annuli  contained in the frontier of a solid torus or
interval bundle component of $\Sigma(M)$, then no annulus in $C$ is homotopic into a
toroidal boundary component of $M$, so $\rho(\pi_1(C))$ is purely hyperbolic. In these cases,
Lemma \ref{solid torus registering}  implies the existence of a registering subgroup for $C$.
Otherwise, $C$ is the frontier of a thickened torus component of $\Sigma(M)$ and 
Lemma \ref{thickened torus registering} guarantees  the existence of a registering subgroup for $C$.
Therefore, we have established:

\begin{prop}
\label{registering exists}
Suppose that $M$ is a compact, orientable, hyperbolizable 3-manifold with incompressible
boundary and $C$ is a characteristic collection of annuli for $M$,
then there exists a $C$-registering subgroup of $\pi_1(M)$.
\end{prop}

\section{Registering subgroups and relative deformation spaces} \label{smaller deformation spaces}

If $H$ is a $C$-registering subgroup for a characteristic collection of annuli $C$,  
then we define $GF(H)$ to be the set of conjugacy classes of 
discrete, faithful, geometrically finite, minimally parabolic representations. (If $H$ is a free group, $GF(H)$ is the space of Schottky representations.)
$GF(H)$ naturally sits inside
$$X_T(H)={\rm Hom}_T(H,\PSL(2,\C))//\PSL(2,\C)$$
where ${\rm Hom}_T(H,\PSL(2,\C))$ denotes the set of representations of
$H$ into $\PSL(2,\C)$ such that
if an element of $H$ lies in a rank two free abelian subgroup of $H$
then its image is either parabolic or the identity.
In turn, $X_T(H)$ is a subvariety of the full character variety $X(H)$.

\begin{lemma}
\label{registering prop disc}
Let $M$ be a compact, orientable, hyperbolizable 3-manifold with non-empty
incompressible boundary.
If $H$ is a registering subgroup for some characteristic collection of annuli,
then
\begin{enumerate}
\item 
$GF(H)$ is an open subset of $X_T(H)$.
\item
If $\{\alpha_n\}$ is a sequence of distinct elements in $\Out(H)$
and $D$ is a compact subset of $GF(H)$, then $\{\alpha_n(D)\}$
exits every compact subset of $X_T(H)$.
\end{enumerate}
\end{lemma}

\begin{proof}
Theorem 10.1 in Marden \cite{marden} implies that $GF(H)$ is an open subset of
$X_T(H)$, which establishes (1).

If (2) fails, there exists 
a sequence $\{\alpha_n\}$ of distinct elements of $\Out(H)$ and a compact
subset  $D$ of $GF(H)$ such that
$\alpha_n(D)$ intersects  a fixed  compact subset of $X_T(H)$ for all $n$. 

We will call an element of $H$ toroidal if it lies in a rank two free abelian subgroup.
Given $\rho\in GF(H)$ and $g\in G$,  let $l_\rho(g)$ denote the translation distance of $\rho(g)$.

Fix, for the moment, an element $\tau\in D$.
Then, given any $P>0$ there exists
finitely many conjugacy classes of  non-toroidal elements $g$ in $H$ such that
\hbox{$l_\tau(g)<P$}. Moreover,  there exists a positive lower
bound on the translation distance $l_\tau(g)$ whenever
$g$ is non-toroidal.
Let $\{h_1,\ldots,h_r\}$ be a generating set for $H$
consisting of non-toroidal elements. If $\{\alpha_n\}$ is a sequence of distinct
elements of $\Out(H)$, then  we may pass to a  subsequence $\{\alpha_j\}$ such that either
\begin{enumerate}
\item
there exists a generator $h_k$ such that 
$l_\tau(\alpha_j^{-1}(h_k))\to\infty$, or
\item
there exist generators $h_i$ and $h_k$
such that the distance between the axes of $\tau(\alpha_j^{-1}(h_i))$ and $\tau(\alpha_j^{-1}(h_k))$ goes to infinity.
\end{enumerate}
In the second case $l_\tau(\alpha_j^{-1}(h_i h_k))\to\infty$. 
Therefore, $\{\alpha_j(\tau)\}$ leaves every compact subset of $X_T(H)$.

Without loss of generality, we may assume that $D$ is contained in a single component of $GF(H)$.
Since all elements  in a component of $GF(H)$ are quasiconformally conjugate 
(see, e.g. \cite[Section 7.3]{CM})
and $D$ is compact,
there exists $L$ such that all the representations in $D$ are $L$-quasiconformally
conjugate.
Therefore, there exists $K>0$ such that all the actions on $\Ht$ are
$K$-bilipschitz conjugate (\cite[Proposition 7.2.6]{CM}). In particular, if $\rho\in D$ and 
$g\in H$, then
$$l_\rho(g)\ge \frac{1}{K}l_\tau(g).$$
It follows that $\{\alpha_j(D)\}$ exits every compact subset of $X_T(H)$, which contradicts
our assumption and completes the proof of (2).
\end{proof}

If $\Sigma$ is an interval bundle component of $\Sigma(M)$,
let $\partial_1\Sigma$ denote the collection of components of $Fr(\Sigma)$ which
are homotopic into toroidal boundary components of $M$.
Let $GF(\Sigma,\partial_1 \Sigma)$  denote the set of conjugacy classes of
discrete, faithful, geometrically finite representations such that the image of a
non-trivial element is parabolic if and only if it is conjugate into $\pi_1(\partial_1\Sigma)$.
$GF(\Sigma,\partial_1\Sigma)$ naturally sits inside
$$X(\Sigma,\partial_1\Sigma)={\rm Hom}(\Sigma,\partial_1\Sigma,\PSL(2,\C))//\PSL(2,\C)$$
where ${\rm Hom}(\Sigma,\partial_1\Sigma,\PSL(2,\C))$
denotes the
representations such that $\rho(g)$ is parabolic or trivial  if $g$ is conjugate into
$\pi_1(\partial_1\Sigma)$.  $X(\Sigma,\partial_1\Sigma)$ is a subvariety of $X(\Sigma)$.

We may use the same argument as in the proof of Lemma \ref{registering prop disc},
replacing non-toroidal elements with elements not conjugate into
$\pi_1(\partial_1\Sigma)$, to establish:

\begin{lemma}
\label{interval bundle prop disc}
Let $M$ be a compact, orientable, hyperbolizable 3-manifold with non-empty
incompressible boundary.
Let $\Sigma$ be an interval bundle component of $\Sigma(M)$, then
\begin{enumerate}
\item
$GF(\Sigma,\partial_1\Sigma)$ is an open subset of $X(\Sigma,\partial_1\Sigma)$.
\item
If $D$ is a compact subset of $GF(\Sigma,\partial_1\Sigma)$ and $\{\alpha_n\}$
is a sequence of distinct elements  of $E(\Sigma,Fr(\Sigma))$, then
$\{\alpha_n(D)\}$ exits every compact subset of $X(\Sigma,\partial_1\Sigma)$.
\end{enumerate}
\end{lemma}

\section{The domains of discontinuity}  \label{domains}

We are now ready to define the domains of discontinuity which occur in the statements of
Theorems \ref{main} and \ref{main absolute}.  We first define $W(M)\subset X_T(M)$.

\begin{defn}
A representation $\rho\in X_T(M)$ lies in $W(M)$ if and only if the following hold:
\begin{enumerate}
\item[(a)]
if $C$ is a characteristic collection of annuli for $M$, then there exists
a \hbox{$C$-registering} subgroup $H$ such that $\rho|_H\in GF(H)$, i.e. $\rho|_H$ is discrete, faithful,
geometrically finite and minimally parabolic, and
\item[(b)]
if $\Sigma$ is an interval bundle component of $\Sigma(M)$ which is not tiny,
then \hbox{$\rho|_{\pi_1(\Sigma)}\in GF(\Sigma,\partial_1\Sigma)$}, i.e. $\rho|_{\pi_1(\Sigma)}$ is discrete, faithful, geometrically finite 
and $\rho|_{\pi_1(\Sigma)}(g)$ is parabolic if and only if $g$ is conjugate to a non-trivial
element of $\pi_1(\partial \Sigma_1)$.
\end{enumerate}
\end{defn}

\begin{prop}
\label{structure of W}
Let $M$ be a compact, orientable, hyperbolizable 3-manifold with non-empty
incompressible boundary which is not an interval bundle.  Then
\begin{enumerate}
\item
$W(M)$ is an $\Out(\pi_1(M))$-invariant open subset of $X_T(M)$.
\item
The interior of $AH(M)$ is  a proper subset of $W(M)$.
\item  
Every minimally parabolic representation in $AH(M)$ lies in $W(M)$. In particular,
$W(M)$ contains a dense subset of $\partial AH(M)$.
\item
$AH(M)\subset W(M)$ if and only if $M$ contains no
primitive essential annuli.
\end{enumerate}
\end{prop}

\begin{proof}{}
We first show that $W(M)$ is open in $X_T(M)$.

Recall that if $H$ is  a registering subgroup for a characteristic
collection of annuli, then \hbox{$r_H:X_T(M)\to X_T(H)$} is continuous.
Since $GF(H)$ is an open subset of $X_T(H)$ (see Lemma \ref{registering prop disc}), 
\hbox{$r_H^{-1}(GF(H))$} is an open subset of 
$X_T(M)$.  Therefore, the set of
representations satisfying condition (a) in the definition of $W(M)$ is open.

If $\Sigma$ is an interval bundle component of $\Sigma(M)$,
then the map \hbox{$r_\Sigma:X_T(M)\to X(\Sigma,\partial_1\Sigma)$} obtained by restriction
is continuous. Since $GF(\Sigma,\partial_1\Sigma)$ is an open subset of $X(\Sigma,\partial_1\Sigma)$
(see Lemma \ref{interval bundle prop disc}),
$r_\Sigma^{-1}(GF(\Sigma,\partial_1\Sigma))$ is an open subset of $X_T(M)$.
It follows that the set of
representations satisfying condition (b) in the definition of $W(M)$ is open.
Therefore, $W(M)$ is open in $X_T(M)$.

Johannson's Classification Theorem implies that every homotopy equivalence $h$ of $M$
is homotopic to one which preserves $\Sigma(M)$ and $M-\Sigma(M)$. In particular, we may assume
that $h$  takes every  interval bundle component of $\Sigma(M)$ to an interval
bundle component of $M$ and takes each characteristic collection of annuli to a characteristic collection of
annuli. Moreover,  if $H$ is a registering subgroup for $C$, we see that $h_*(H)$ is a registering subgroup
for $h(C)$. Since every outer automorphism of $\pi_1(M)$ is realized by a homotopy equivalence,
one easily verifies that $W(M)$ is invariant under $\Out(\pi_1(M))$. This completes the proof of (1).

Since all representations in ${\rm int}(AH(M))$ are minimally parabolic,
(2) follows from (3). We now turn to the proof of (3).

Suppose that $\rho\in AH(M)$ is minimally parabolic. 
Lemmas \ref{solid torus registering} and \ref{thickened torus registering} imply
that if $C$ is a characteristic collection of annuli, then there is a \hbox{$C$-registering}
subgroup $H$ such that $\rho|_H$ is discrete, faithful, geometrically finite and minimally parabolic.
Therefore,  $\rho$ satisfies condition (a) in the definition of $W(M)$.

Now suppose that $\Sigma$ is an interval bundle component of $\Sigma(M)$ which is
not tiny.  Let $M_\rho$ be a relative compact core for $N^0_\rho$ and let $h:M\to M_\rho$ be a homotopy
equivalence in the homotopy class of $\rho$. Johannson's Classification Theorem implies
that we may assume that $h(\Sigma)=\Sigma_\rho$ is an interval bundle component of
$\Sigma(M_\rho)$ and that $h$ restricts to a homeomorphism from $Fr(\Sigma)$ to $Fr(\Sigma_\rho)$.
Let $\partial_1\Sigma_\rho=h(\partial_1\Sigma)$.
Since $\rho$ is minimally parabolic, $r_\Sigma(\rho)(g)$ is parabolic if and only if $g$ is conjugate
to a non-trivial element of $\pi_1(\partial_1\Sigma)$. 

The interval bundle $\Sigma_\rho$ lifts to a compact core 
for the cover $N_\Sigma$ of $N_\rho$ associated to $\rho(\pi_1(\Sigma))=\pi_1(\Sigma_\rho)$.
However, the lift need not be a relative compact core, since it need not intersect every component of
$N^0_\Sigma$ in an incompressible annulus.
(Here we choose the invariant system of horoballs for $\rho(\pi_1(\Sigma))$
to be a subset of the precisely invariant system of horoballs for $\rho(\pi_1(M))$, so the covering map
from $N_\Sigma$ to $N_\rho$ restricts to a covering map from $\partial N^0_\Sigma$ to its image in
$\partial N^0_\rho$.) In order to extend $\Sigma_\rho$ to a submanifold which does lift
to a relative compact core, we construct a submanifold $Y_\rho$ of $M_\rho$ which is homeomorphic to
$\partial_1\Sigma_\rho\times [0,1]$  by a homeomorphism identifying $\partial_1\Sigma_\rho$
with $\partial_1\Sigma_\rho\times \{0\}$, so that $Y_\rho\cap \Sigma_\rho=\partial_1\Sigma_\rho$ and
$Y_\rho\cap \partial N^0_\rho$ is a collection of incompressible annuli which is identified with
$\partial_1\Sigma_\rho\times \{1\}$.
If we let
$$\Sigma_\rho^+=\Sigma_\rho\cup Y_\rho$$
then $\Sigma_\rho^+$ does lift to a relative compact core for $\partial N_\Sigma^0$.
Moreover, the lift of $\Sigma_\rho^+$ intersects $N^0_\Sigma$ exactly in
the lift of $\Sigma_\rho^+\cap \partial N^0_\rho$.
The ends of $N_\Sigma^0$ are in one-to-one correspondence with the components of
$\partial\Sigma_\rho^+-(\partial \Sigma_\rho^+\cap \partial N^0_\rho)$, each of which is homotopic to a component of
$\partial \Sigma_\rho-\partial_1\Sigma$. In particular, $N^0_\Sigma$ has one or two ends.

If the manifold $N^0_\Sigma$ has only one end, then the covering map
\hbox{$N_{\Sigma}\to N_\rho$} is infinite-to-one on this end. The Covering
Theorem  implies that  the single end of $N^0_{\Sigma}$ is geometrically finite, so
$N_\Sigma$ is geometrically finite.
If the manifold $N^0_\Sigma$ has two ends, then
$Fr(\Sigma_\rho)=\partial_1\Sigma_\rho$, so
each component of 
$\partial \Sigma_\rho-\partial_1\Sigma$
is identified with a proper subsurface of a component of $\partial M_\rho.$
Again \hbox{$N_{\Sigma}\to N_\rho$} 
is infinite-to-one on each end of $N^0_\Sigma$
and the Covering Theorem may be used to show that $N_\Sigma$ is geometrically finite.
Thus, in all cases, $\rho|_{\pi_1(\Sigma)}\in GF(\Sigma,\partial_1\Sigma)$, so $\rho$ satisfies condition (b) in
the definition of $W(M)$. Therefore, minimally parabolic representations in $AH(M)$ lie in $W(M)$.

Since minimally parabolic representations are dense in the boundary
of $AH(M)$,
$W(M)$ contains a dense subset of $\partial AH(M)$.
(The density of minimally parabolic representations in $\partial AH(M)$ follows from
Lemma 4.2 in \cite{canary-hersonsky}, which shows that minimally parabolic
representations are dense in the boundary of any component of ${\rm int}(AH(M))$
and the Density Theorem, see Brock-Canary-Minsky \cite{BCM}, Bromberg-Souto
\cite{bromberg-souto}, Namazi-Souto \cite{namazi-souto}
or Ohshika \cite{ohshika-density}, which asserts that $AH(M)$ is the closure of its interior.)
This  completes the proof of (3). 

Suppose that
$M$ contains no primitive essential annuli. Then $\Sigma(M)$ contains 
no interval bundle components which are not tiny, since otherwise a non-peripheral
essential annulus in the interval bundle would be a primitive essential annulus
(see \cite[Lemma 7.3]{canary-storm2}).
Similarly, every component
of the frontier of a tiny interval bundle component is isotopic into a 
solid torus or thickened
torus component of $\Sigma(M),$ since otherwise it would be a primitive essential
annulus (see Johannson \cite[Lemma 32.1]{johannson}).
Therefore, every characteristic collection of annuli is the frontier of either a
solid torus or thickened torus component of $\Sigma(M)$.
Moreover, the core curve of
each solid torus component $V$  of $\Sigma(M)$ is non-peripheral, since otherwise
its frontier annuli would be primitive essential annuli (again see
Johannson \cite[Lemma 32.1]{johannson}).
Therefore,
just as in the proof of Lemma 8.1 in \cite{canary-storm2}, $\rho(\pi_1(V))$
is purely hyperbolic if $\rho\in AH(M)$ and $V$ is a solid torus component of $\Sigma(M)$.
Therefore, if $\rho\in AH(M)$ and $C$ is any characteristic collection of annuli for $M$,
then Lemma \ref{solid torus registering} or \ref{thickened torus registering}
guarantees that there exists a
\hbox{$C$-registering} subgroup $H$  such that $\rho|_H$ is discrete, faithful,
geometrically finite and minimally parabolic. Since every interval bundle is tiny,
it follows that $AH(M)\subset W(M)$.

On the other hand, if $M$ contains a primitive essential annulus $A$, then
there exists $\rho\in AH(M)$ such that $\rho(\pi_1(A))$ is purely parabolic (see Ohshika
\cite{ohshika-parabolic}). Since $A$ is either isotopic to a component of a characteristic collection of annuli
or isotopic into an interval bundle component of $\Sigma(M)$,
$\rho$ does not lie in $W(M)$. In particular, $AH(M)$ is not a
subset of $W(M)$. Therefore, $AH(M)\subset W(M)$ if and only if $M$
contains no primitive essential annuli.

\end{proof}

If $M$ does not contain an essential annulus which intersects a toroidal boundary
component of $M$, then we define $\hat W(M)\subset X(M)$.  Notice that $M$ contains an essential annulus with
one boundary component contained in a toroidal boundary
component of $M$ if and only if
$\Sigma(M)$ has a thickened torus component.

\begin{defn}
We say that
$\rho\in X(M)$ lies in $\hat W(M)$ if and only if the following hold:
\begin{enumerate}
\item[(a)]
if $C$ is a characteristic collection of annuli for $M$, then there exists
a \hbox{$C$-registering} subgroup $H$ such that $\rho|_H\in GF(H)$, and
\item[(b)]
if $\Sigma$ is an interval bundle component of $\Sigma(M)$ with base surface $F$, which is not tiny,
then $\rho|_{\pi_1(\Sigma)}\in GF(\Sigma,\emptyset)$, i.e. $\rho|_{\pi_1(\Sigma)}$
is discrete, faithful, geometrically finite and purely hyperbolic.
\end{enumerate}
\end{defn}

We obtain the following analogue of Proposition \ref{structure of W} whenever $\hat W(M)$ is
defined.

\begin{prop}
\label{structure of W hat}
Let $M$ be a compact, orientable, hyperbolizable 3-manifold with non-empty
incompressible boundary which is not an interval bundle and so that $M$ contains
no essential annuli which intersects a toroidal boundary component.
Then
\begin{enumerate}
\item
$\hat W(M)$ is an $\Out(\pi_1(M))$-invariant open subset of $X(M)$.
\item 
$\hat W(M)\cap X_T(M)=W(M)$.
\item
$AH(M)\subset \hat W(M)$ if and only if $M$ contains no
primitive essential annuli.
\end{enumerate}
\end{prop}
 
\medskip\noindent
{\em Sketch of proof:}  
Since $\Sigma(M)$ does not contain any thickened torus components, every characteristic
collection of annuli is either the frontier of a solid torus component of $\Sigma(M)$ or a component
of the frontier of an interval bundle component of $\Sigma(M)$. Moreover, every registering subgroup
$H$ is a free group and $\partial_1\Sigma$ is empty for every interval bundle component of $\Sigma(M)$.
The proof of (1) mimics the proof of  Proposition \ref{structure of W}. 
If $H$ is a registering subgroup for some characteristic collection of annuli, then we can define
$r_H:X(M)\to X(H)$, and $r_H^{-1}(GF(H))$ is an open subset of $X(M)$. 
In the case that $\Sigma$ is an interval bundle component of $\Sigma(M)$, we define
\hbox{$r_\Sigma:X(M)\to X(\Sigma)$} and $GF(\Sigma,\emptyset)$
is an open subset of $X(\Sigma)$, so $\rho_\Sigma^{-1}(GF(\Sigma,\emptyset))$ is open in
$X(M)$.  Therefore, as in the proof of property (1) in Proposition \ref{structure of W}, $W(M)$ is
an open subset of $X(M)$.  The 
$\Out(\pi_1(M))$-invariance of $\hat W(M)$ follows from Johannson's Classification Theorem,
much as in the proof of Proposition \ref{structure of W}.
 
Property (2) follows immediately from the definitions of $W(M)$ and $\hat W(M)$ and the restrictions
on the characteristic submanifold of $M$ discussed in the previous paragraph.  Property (3) follows from
property (2) and part (4) of Theorem \ref{structure of W}.
\eproof

\noindent
{\bf Remark:} If one, more generally, allowed $\rho|_H$ and $\rho|_{\pi_1(\Sigma)}$ to be
primitive-stable (see Minsky \cite{minsky-primitive}) in the definition of $\hat W (M)$,
then $\hat W(M)$ would agree with the domain of discontinuity obtained in 
\cite{canary-storm2} in the case that $M$ has no toroidal boundary components.

\section{Proof of Main Theorem} \label{prop disc section}

We are now prepared to complete the proof of our main theorem, which we recall below:

\bigskip
\noindent
{\bf Theorem \ref{main}.} {\em Let $M$ be a compact, orientable, hyperbolizable $3$-manifold with nonempty incompressible boundary, which is not an
interval bundle. Then there exists an open \hbox{${\rm Out}(\pi_1(M))$}-invariant subset 
$W(M)$ of  $X_T(M)$ such that \hbox{${\rm Out}(\pi_1(M))$}  acts properly discontinuously on $W(M)$, \hbox{${\rm int}(AH(M))$} is a proper subset of $W(M)$, and
$W(M)$ intersects $ \partial AH(M)$.}

\medskip

Theorem \ref{main} follows  immediately from Proposition \ref{structure of W}, which gives the key properties
of $W(M)$, and the following
proposition which establishes the proper discontinuity of the action of $\Out(\pi_1(M))$ on $W(M)$.

\begin{prop}
\label{prop disc on W}
If $M$ is a compact, orientable, hyperbolizable  3-manifold with non-empty
incompressible boundary which is not an interval bundle, then
$\Out(\pi_1(M))$ acts properly discontinuously on $W(M)$.
\end{prop}

\begin{proof}
Since $\Out_0(\pi_1(M))$ has finite index in $\Out(\pi_1(M))$,  it suffices
to prove that $\Out_0(\pi_1(M))$ acts properly discontinuously on $W(M)$.

Suppose that $\Out_0(\pi_1(M))$ does not act properly discontinuously on $W(M)$.
Then there exists a compact subset $R$ of $W(M)$ and a sequence
$\{\alpha_n\}$ of distinct elements in $\Out_0(\pi_1(M))$ such that
$\alpha_n(R)\cap R$ is non-empty for all $n$.
We may pass to a subsequence so that either 
\begin{enumerate}
\item
$\{\Phi(\alpha_n)\}$ is a sequence of distinct elements of
$\left(\oplus_i D(T_i)\right)\oplus\left(\oplus_j E(\Sigma_j,Fr(\Sigma_j)\right)$, or
\item
$\{\Phi(\alpha_n)\}$ is a constant sequence.
\end{enumerate}

In case (1) we may pass to a further subsequence, still called
$\{\alpha_n\}$, so that either
\begin{enumerate}
\item[(a)]
there exists an interval bundle component $\Sigma$ of $\Sigma(M)$ so
that 
$\{p_\Sigma(\Phi(\alpha_n)\}$ is a sequence of distinct elements of $E(\Sigma,Fr(\Sigma))$
where 
$$p_{\Sigma}:\left(\oplus_i D(T_i)\right)\oplus\left(\oplus_j E(\Sigma_j,Fr(\Sigma_j)\right)\to E(\Sigma, Fr(\Sigma)),$$
is the obvious projection map onto $E(\Sigma,Fr(\Sigma)$, or
\item[(b)]
there exists a thickened torus component $T$ of $\Sigma(M)$
so
that 
$\{p_T(\Phi(\alpha_n)\}$ is a sequence of distinct elements of $D(T)$
where 
$$p_T:\left(\oplus_i D(T_i)\right)\oplus\left(\oplus_j E(\Sigma_j,Fr(\Sigma_j)\right)\to D(T)$$
is the obvious projection map onto $D(T)$.
\end{enumerate}

In case (1a), $r_\Sigma(R)$ is a compact subset of $GF(\Sigma,\partial_1\Sigma)$ and
$\{p_\Sigma(\Phi(\alpha_n)\}$ is a sequence of distinct elements of
$E(\Sigma,Fr(\Sigma))$. Recall, see Lemma \ref{interval bundle E}, that
$E(\Sigma,Fr(\Sigma))$ is identified with a subgroup of $\Out(\pi_1(\Sigma))$.
Notice that, by construction,
$$r_\Sigma(\alpha(\rho))=p_\Sigma(\Phi(\alpha))(r_\Sigma(\rho)).$$
or stated differently,
$$r_\Sigma(\rho \circ \alpha^{-1}) = r_\Sigma(\rho) \circ p_\Sigma(\Phi(\alpha))^{-1}.$$
for all
$\rho\in W(M)$ and all $\alpha\in \Out_0(\pi_1(M))$. Therefore,
$p_\Sigma(\Psi(\alpha_n))(r_\Sigma(R))\cap r_\Sigma(R)$ is non-empty for all $n$.
Since $r_\Sigma(R)$ is a compact subset of $GF(\Sigma,\partial_1\Sigma)$,
this contradicts the proper discontinuity of the action of
$E(\Sigma,Fr(\Sigma))$ on $GF(\Sigma,\partial_1\Sigma)$, see Lemma  \ref{interval bundle prop disc}.
This contradiction rules out case (1a).  

In case (1b), notice that
if $\rho\in W(M)$, then $\rho|_{\pi_1(T)}:\pi_1(T)\to\PSL(2,\C)$ is discrete
and faithful. Therefore, there exists a continuous restriction map \hbox{$r_T:W(M)\to AH(\BZ^2)$}, 
where $AH(\BZ^2)$ is the space
of conjugacy classes of discrete faithful representations from $\pi_1(T)$ into
$\PSL(2,\C)$. There is a natural identification of $D(T)$ with a subgroup of 
$\Out(\pi_1(T))$.
Since all the elements of ${\rm ker}(p_T\circ \Phi)$ act trivially on $\pi_1(T)$,
$$r_T(\alpha(\rho)) = p_T(\Phi(\alpha))(r_T(\rho)).$$
for all $\rho\in W(M)$ and $\alpha\in \Out_0(\pi_1(M))$.
Therefore, $p_T(\Phi(\alpha_n))(r_T(R))\cap r_T(R)$ is non-empty for all $n$.
It is easy to check that $\Out(\pi_1(T))$ acts properly discontinuously
on $AH(\BZ^2)$. (One may identify $AH(\BZ^2)$ with $\mathbb{C} \backslash \mathbb{R}$ and the action of
$\Out(\pi_1(T))$ is identified with the action of ${\rm GL}(2,{\mathbb{Z}})$ as a group of conformal and anti-conformal
automorphisms of $\mathbb{C} \backslash \mathbb{R}$,
 which  is well-known to act properly discontinuously on 
$\mathbb{C} \backslash \mathbb{R}$.)
So, we have again obtained a contradiction and case (1b) cannot occur.

In case (2), there exists $\gamma\in \Out_0(\pi_1(M))$ and a sequence $\beta_n\in \ker(\Phi)$
such that \hbox{$\alpha_n=  \beta_n\circ\gamma$} for all $n$. 
Since $\gamma$ induces a homeomorphism of $X_T(M)$ which preserves $W(M)$,
$\gamma(R)$ is a compact subset of $W(M)$ and
$\beta_n(\gamma(R))\cap R$ is non-empty for all $n$.
Recall that $\ker(\Phi)=\oplus \hat K(C_j)$. Therefore, after passing to a further subsequence,
we may find a characteristic collection of annuli $C$ so that $q_C(\beta_n)$ is
a sequence of distinct elements of $\hat K(C)$ (where $q_C$ is the projection of
$\ker(\Phi)$ onto $\hat K(C)$).
Since $X_T(M)$ is locally compact, 
for each $x\in W(M)$, there exists an open neighborhood $U_x$
of $x$ and a \hbox{$C$-registering} subgroup $H_x$ such that the closure
$\bar U_x$ is a compact subset of $W(M)$ and $r_{H_x}(\bar U_x)\subset GF(H_x)$.
Since $\gamma(R)$ is compact,
there exists a finite collection of points \hbox{$\{ x_1,\ldots, x_r\}$}
such that
\hbox{$\gamma(R)\subset U_{x_1}\cup \cdots\cup U_{x_r}$}. Therefore, again passing
to subsequence if necessary, there must exist
$x_i$ such that \hbox{$\beta_n(U_{x_i})\cap R$} is non-empty for all $n$. Let
$U=U_{x_i}$ and $H=H_{x_i}$.
Lemma \ref{characteristic} implies that 
\hbox{$\{s_{H}(q_C(\beta_n))\}$} is a sequence of distinct elements of \hbox{${\rm Out}(H)$}
and that
$$s_{H}(q_C(\beta_n))(r_{H}(\bar U))=r_{H}(\beta_n(\bar U)).$$
Lemma \ref{registering prop disc} then
implies that  
$$\{s_{H}(q_C(\beta_n))(r_{H}(\bar U))\}=\{ r_{H}(\beta_n(\bar U))\}$$
exits every compact
subset of \hbox{$X_T(H)$}. Therefore, \hbox{$\{\beta_n(U)\}$} exits every compact subset of
$X_T(M)$ which is  again a  contradiction.
Therefore, case (2) cannot occur and we have completed the proof.
\end{proof}

Corollary \ref{maincor}, which we restate here, follows readily from Theorem \ref{main},
Proposition \ref{structure of W} and Theorem 1.2 from \cite{canary-storm2}.

\medskip\noindent
{\bf Corollary \ref{maincor}:} {\em
If  $M$ is a compact, orientable, hyperbolizable  3-manifold with incompressible boundary and non-abelian fundamental group,
then $\Out(\pi_1(M))$ acts properly discontinuously on an open
\hbox{${\rm Out}(\pi_1(M))$}-invariant
neighborhood of $AH(M)$ in $X_T(M)$
if and only if $M$ contains no primitive essential
annuli.}

\medskip\noindent
{\em Proof of Corollary \ref{maincor}:}
Proposition \ref{prop disc on W} shows that $\Out(\pi_1(M))$ acts properly discontinuously on
$W(M)$, and Proposition \ref{structure of W} implies that $W(M)$ is an open neighborhood of $AH(M)$ when $M$ contains no primitive essential annuli. 

If $M$ contains a primitive essential annulus, then Theorem 1.2 of \cite{canary-storm2} asserts
that $\Out(\pi_1(M))$ does not act properly discontinuously on $AH(M)$, so $\Out(\pi_1(M))$ cannot act properly
discontinuously on any \hbox{${\rm Out}(\pi_1(M))$}-invariant
neighborhood of $AH(M)$ in $X_T(M)$.
\eproof

\noindent
{\bf Remark:} One could also prove that $\Out(\pi_1(M))$ cannot act properly discontinuously
on an open neighborhood of $AH(M)$ when $M$ contains a primitive essential annulus using the technique of Lemma 15 in Lee \cite{lee}.

\section{Dynamics in the absolute character variety} \label{final cor section}

In this section, we study the action of $\Out(\pi_1(M))$ on the full character variety $X(M)$.
We begin by showing that $\Out(\pi_1(M))$ acts properly discontinuously on $\hat W(M)$,
which is a nearly immediate generalization of
Proposition \ref{prop disc on W}.

\begin{prop}
\label{prop disc on hat W}
If $M$ is a compact, orientable, hyperbolizable  3-manifold with non-empty
incompressible boundary which is not an interval bundle,  and no essential annulus
in $M$ has a boundary component contained in a toroidal boundary component of $M$, then
$\Out(\pi_1(M))$ acts properly discontinuously on $\hat W(M)$.
\end{prop}

\medskip\noindent
{\em Sketch of proof:}  Again, it suffices to prove that $\Out_0(\pi_1(M))$ acts properly
discontinuously on $\hat W(M)$. 
If  $\Out_0(\pi_1(M))$ does not act properly discontinuously on $\hat W(M)$,
then there exists a compact subset $R$ of $\hat W(M)$ and a sequence
$\{\alpha_n\}$ of distinct elements in $\Out_0(\pi_1(M))$ such that
$\alpha_n(R)\cap R$ is non-empty for all $n$. We may again pass to a subsequence so
that either (1) $\{\Phi(\alpha_n)\}$ is a sequence of distinct elements or (2) $\{\Phi(\alpha_n\}$
is a constant sequence. Since $\Sigma(M)$ contains no thickened torus components, in case (1)
we can assume that there exists an interval bundle component $\Sigma$ of $\Sigma(M)$ such
that $\{p_\Sigma(\alpha_n)\}$ is a sequence of distinct elements of $E(\Sigma,Fr(\Sigma))$.
We then proceed, exactly as in the consideration of cases (1)(a) and (2) in the proof of
Proposition \ref{prop disc on W}, to obtain a contradiction.
\eproof

Propositions \ref{structure of W hat} and \ref{prop disc on hat W} immediately
imply Theorem \ref{main absolute}.

\medskip\noindent
{\bf Theorem \ref{main absolute}:} {\em
Let $M$ be a compact, orientable, hyperbolizable $3$-manifold with nonempty incompressible boundary, which is not an interval bundle. If $M$ does not contain an essential annulus with
one boundary component contained in a toroidal boundary component of $M$,
then there exists an open \hbox{${\rm Out}(\pi_1(M))$}-invariant subset 
$\hat W(M)$ of  $X_T(M)$ such that \hbox{${\rm Out}(\pi_1(M))$}  acts properly discontinuously on $\hat W(M)$ and
$$W(M)=\hat W(M)\cap X_T(M).$$
In particular, $\hat W(M)$ intersects $\partial AH(M)$.}

\bigskip

We next adapt the proof of Lemma 15 in Lee \cite{lee} to establish Proposition
\ref{bad annuli}:

\medskip\noindent
{\bf Proposition \ref{bad annuli}:} {\em
Let $M$ be a compact, orientable, hyperbolizable $3$-manifold with nonempty incompressible boundary and non-abelian fundamental group.
If $M$ contains an essential annulus with
one boundary component contained in a toroidal boundary component, then every point in $AH(M)$ is a limit
of representations in $X(M)$ which are fixed points of infinite order elements of $\Out(\pi_1(M))$.}

\medskip
\noindent
\begin{proof}
Let $A$ be an essential annulus in $M$ with one boundary component contained in
a toroidal boundary component of $M$ and let $a$ be the core curve of $A$.

First suppose that  $\rho \in {\rm int}(AH(M))$. Theorem 5.7 in Bromberg
\cite{bromberg-rigidity} implies that there exists a 
neighborhood $U$ of $\rho\in X(M)$ and an open holomorphic
map \hbox{$Tr_a:U\to \C$} such that if $\rho'\in U$, then
the trace of $\rho'(a)$ is given by $\pm Tr_a(\rho')$. (Recall that the trace of a representation into
$\PSL(2,\C)$ is only well-defined up to sign.)
Therefore, there
exists a sequence $\{\rho_n\}\subset X(M)$ such that $\{\rho_n\}$ converges to $\rho$ and
$\rho_n(a)^n=Id$ for all large enough $n$. (Simply choose a sequence of representations
$\{\rho_n\}$ converging to $\rho$ such that $Tr_a(\rho_n)=\pm 2\cosh(\frac{\pi}{n})$.)
For each $n$, $\rho_n$ is fixed by
the infinite order element \hbox{$(D_A)_*^{n}\in\Out(\pi_1(M))$} 
where $D_A$ is the Dehn twist about $A$. Therefore, $\rho$ is a limit of
fixed points of infinite order elements of $\Out(\pi_1(M))$.

The Density Theorem (\cite{BCM,bromberg-souto,namazi-souto,ohshika-density}) assures
us that $AH(M)$ is the closure of its interior, so, by diagonalization,
every representation in $AH(M)$  is also a limit of
fixed points of infinite order elements of $\Out(\pi_1(M))$.
\end{proof}

One may combine Proposition \ref{bad annuli} with
Proposition \ref{structure of W hat}, Theorem \ref{main absolute} and Theorem 1.2
from \cite{canary-storm2} to prove Corollary \ref{maincor absolute}.

\medskip\noindent
{\bf Corollary \ref{maincor absolute}:}  {\em
If  $M$ is a compact, orientable, hyperbolizable   3-manifold with incompressible boundary and non-abelian fundamental group, then  $\Out(\pi_1(M))$ acts properly discontinuously on an open,
\hbox{${\rm Out}(\pi_1(M))$}-invariant
neighborhood of $AH(M)$ in $X(M)$  if and only if $M$  does not contain a primitive essential annulus or
an essential annulus with one boundary component contained in a toroidal boundary
component of $M$.}

\begin{proof}
Suppose that $M$  does not contain a primitive essential annulus or
an essential annulus with one boundary component contained in a toroidal boundary
component of $M$.
Theorem \ref{main absolute} implies that $\Out(\pi_1(M))$ acts properly discontinuously on
$W(M)$, while Proposition \ref{structure of W hat} implies that $AH(M)\subset W(M)$ and
that  $W(M)$ is open in $X(M)$. 
Therefore, $\Out(\pi_1(M))$ acts properly discontinuously on an open neighborhood
of $AH(M)$ in $X(M)$. 

On the other hand, if $M$ contains a primitive essential annulus, then
$\Out(\pi_1(M))$ does not act properly discontinuously on $AH(M)$, by
Theorem 1.2 of \cite{canary-storm2}, so it cannot act properly discontinuously
on an open $\Out(\pi_1(M))$-invariant neighborhood of $AH(M)$. 
If $M$ contains an essential annulus
with a boundary component contained in a toroidal boundary component of $M$,
Proposition \ref{bad annuli}
shows that no point in $AH(M)$ can be contained in a domain of discontinuity
for the action of  $\Out(\pi_1(M))$  on $W(M)$. 
The consideration of these two cases completes the proof of Corollary \ref{maincor absolute}.
\end{proof}

 \end{document}